\documentclass[10pt]{article}

\usepackage{geometry}

\usepackage{parskip}

\usepackage{amsthm}
\usepackage{amsmath}
\usepackage{amsfonts}
\usepackage{amssymb}

\allowdisplaybreaks

\usepackage{algorithm}
\usepackage{algorithmic}
\newcommand{\algorithmicbreak}{\textbf{break}}
\newcommand{\BREAK}{\STATE \algorithmicbreak}

\usepackage{graphicx}
\graphicspath{{figures/}}
\usepackage{hyperref}
\usepackage{color}
\newtheorem{definition}{Definition}[section]

\newtheorem{remark}[definition]{Remark}

\newtheorem{lemma}[definition]{Lemma}

\newcommand{\R}{\mathbb{R}}
\newcommand{\G}{\mathcal{G}}
\newcommand{\C}{\mathcal{C}}

\newcommand{\vecsym}[1]{\boldsymbol{#1}}
\renewcommand{\vec}[1]{\mathbf{#1}}

\DeclareMathOperator{\Tr}{Tr}

\title{Kinetic Methods for Inverse Problems}
\date{\today}
\author{
	Michael Herty  and Giuseppe Visconti \medskip\\
	{\small\it Institut f\"{u}r Geometrie und Praktische Mathematik (IGPM)} \\
	{\small\it RWTH Aachen University} \\
	{\small\it Templergraben 55, 52062 Aachen, Germany
	}
}

\begin{document}

\maketitle

\begin{abstract}
	The Ensemble Kalman Filter method can be used as an iterative numerical scheme for parameter identification or	nonlinear filtering problems. We study the limit of infinitely large ensemble size and derive the corresponding mean-field limit of the ensemble method. The solution of the inverse problem is provided by the expected value of the distribution of the ensembles and the kinetic equation allows, in simple cases, to analyze stability of these solutions. Further, we present a slight but stable modification of the method which leads to a Fokker-Planck-type kinetic equation. The kinetic methods proposed here are able to solve the problem with a reduced computational complexity in the limit of a large ensemble size. We illustrate the properties and the ability of the kinetic model to provide solution to inverse problems by using examples from the literature.
\end{abstract}



\section{  Introduction  }
We are concerned with the following abstract inverse problem or parameter identification problem
\begin{equation} \label{eq:noisyProb}
	\vec{y} = \G(\vec{u}) + \vecsym{\eta}
\end{equation}
where $\G:X \to Y$ is the (possible nonlinear) forward operator between finite dimensional Hilbert spaces $X=\R^d$ and $Y=\R^K$, with $d,K\in\mathbb{N}$, $\vec{u}\in X$ is the control, $\vec{y}\in Y$ is the observation and $\vecsym{\eta}$ is  observational noise. Given noisy measurements or observations $\vec{y}$ and the known mathematical model $\G$, we are interested in finding the corresponding control $\vec{u}$. Typically, the observational noise $\vecsym{\eta}$ is not explicitly known but only information on its distribution is available. Inverse problems, in particular in view of a possible
ill-posedness, have been discussed in vast amount of literature and we refer to 
~\cite{EnglHankeNeubauer1996} for an introduction and further references. In the following
we will investigate a particular numerical method for solving problem~\eqref{eq:noisyProb}, namely, the Ensemble Kalman Filter (EnKF). While this method has already been introduced more than ten years ago~\cite{Evensen1994}, recent theoretical 
progress~\cite{schillingsstuart2017} is the starting point of this work.  

As in \cite{schillingsstuart2017} we aim to solve
the inverse problem by minimizing the  least squares functional
\begin{equation} \label{eq:leastSqFnc}
	\Phi(\vec{u},\vec{y}) := \frac12 \left\| \vecsym{\Gamma}^{\frac12} (\vec{y} - \G(\vec{u})) \right\|^2
\end{equation}
where $\vecsym{\Gamma}^{-1}$ normalizes the so-called model-data misfit. This is defined as the covariance of the noise $\vecsym{\eta}$. Note that there is no regularization of the control $\vec{u}$ in the minimization problem of~\eqref{eq:leastSqFnc}. See e.g.~\cite{BianchiBucciniDonatelliSerra2015,Groetsch1984,Hansen1998,KlannRamlau2008} for examples of Tikhonov and other regularization technique. 

We 
briefly recall a Bayesian inversion formulation for problem~\eqref{eq:noisyProb}. 
Following~\cite{dashtistuart2017,Stuart2010} a solution to the  inverse problem is obtained by treating the unknown control $\vec{u}$, the data $\vec{y}$ and the noise $\vecsym{\eta}$  as random variables. Then, the conditional probability measure of the control $\vec{u}$ given the observation $\vec{y}$, called posterior measure, is computed via Bayes Theorem. Typically, there is an interest in moments of the posterior, e.g. choosing the point of maximal probability (MAP estimator). For further details concerning Bayesian inversion, e.g. the modeling of the unknown prior distributions and other choices of estimators, see~\cite{Berger1985,BurgerLucka2014,dashtistuart2017,ernstetal2015} and references therein. 

Before finally stating the aim of this work, we briefly recall some references on the EnKF method without aiming to give a complete list. Iterative filtering methods have also been successfully applied to inverse problems since many years. A particular successful method has been originally proposed in~\cite{Kalman1960} to estimate state variables, parameters, etc. of stochastic dynamical systems. This method has been extended to the EnKF in~\cite{Evensen1994}. The EnKF sequentially updates each member of an ensemble of random elements in the space $X$ by means of the Kalman update formula, using the knowledge of the model $\G$ and of given observational data $\vec{y}$. It is important to note that {\em no} information on the derivative of $\G$ is required. The EnKF provides satisfactory results even when used with a small number of ensembles, as proved by the accuracy analysis in~\cite{MajdaTong2018}. Some examples in mathematical literature  of the  application of the filtering method to inverse problems are given in the incomplete list~\cite{Oliveretal,schillingsetal2018,SchillingsPreprint,iglesias2015,iglesiaslawstuart2013,schillingsstuart2017,schillingsstuart2018}. In particular, we refer to the following books~\cite{Evensen2009book,OliverReynoldsLiu2008}. Our starting point is~\cite{schillingsstuart2017} where the continuous time limit of the EnKF has been studied as a regularization technique for minimization of the least
squares functional~\eqref{eq:leastSqFnc} with a finite ensemble size. Recently, further study has been conducted in this direction~\cite{ChadaStuartTong2019,LangeStannat2019}. We also note that the EnKF can also formally be derived within the Bayesian framework~\cite{ernstetal2015,iglesiaslawstuart2013b,kwiatkowskimandel2015,LawStuart2012,leglandmonbettran2009}.

In the cited references the ensemble size was fixed and, due to the possible associated high computational cost, limited to a small number of ensembles. The analysis of the method for a large ensemble size limit has been investigated in~\cite{DelMoralKurtzmannTugaut2017,DelMoralTugaut2018,ernstetal2015,lawtembinetempone2016}. However, to the best of our knowledge, an evolution equation for the probability distribution of the unknown control has not been derived. We aim to provide a continuous representation of the EnKF method that also 
holds in the limit of infinitely many ensembles. We do believe that the derivation of the kinetic equation leads to insight to the method that might not be easy to obtain
otherwise.
The main advantage of the derivation of a mean-field equation is twofold. First, it formally allows to deal with the case of infinitely many ensembles
and in this regime numerical simulations show a better reconstruction of an estimator of the unknown control than considering small ensemble sizes. Second, it allows to study stability, at least in the simple case of a one-dimensional control, and it suggests a modification of the method which results in improved stability of the corresponding Fokker-Planck equation.

We proceed as follows: We start from the continuous time limit formulation of the EnKF derived in~\cite{schillingsstuart2017} and interpret it as an interacting particle system. Then, we study the mean-field limit for large ensemble sizes. From a mathematical point of view, this technique has been widely used to reduce the computational complexity and to analyze interacting particle models, e.g. in socio-economic dynamics or gas dynamics~\cite{CarrilloFornasierToscaniVecil2010,CarrilloPareschiZanella2019,CristianiPiccoliTosin2014,hatadmor2008,HertyRinghofer2011,PareschiToscaniBOOK,Toscani2006,TrimbornPareschiFrank}. The kinetic equation evolves in time the probability distribution of the control and the solution to the inverse problem is  shown to be the mean of this distribution. 
We analyze linear stability of the EnKF. Further, we present suitable modifications of the method based on the kinetic formulation in order to improve the stability pattern. The kinetic model guarantees a computational gain in the numerical simulations using a Monte Carlo approach similar to~\cite{AlbiPareschi2013,BabovskyNeunzert1986,FornasierEtAl2011,Lemou1998,MouhotPareschi2006,PareschiRusso1999,PareschiToscaniBOOK}. 


\section{From the Ensemble Kalman Filter to the gradient descent equation} \label{sec:enkf}

The Ensemble Kalman Filter (EnKF) has been introduced~\cite{Evensen1994} as a discrete time method to estimate state variables, parameters, etc. of stochastic dynamical systems. The estimations are based on system dynamics and measurement data
that are possibly perturbed by known noise. The EnKF is a generalization and improved version of the classical Kalman Filter method~\cite{Kalman1960}. In the following, we briefly review the definition of the EnKF which is based on a sequential update of an ensemble of states and parameters. Then we recover the continuous time limit equation derived in the recent work~\cite{schillingsstuart2017}. This will be the starting point to introduce and compute in the next sections a mean-field limit 
for infinitely many ensembles. The arising kinetic partial differential equation allows subsequent analysis on the nature
of the method.

As in~\cite{schillingsstuart2017} we  consider a control $\vec{u}\in\R^d$, a given state $\vec{y} \in \R^K$ coupled by the system
dynamic $\mathcal{G}$ as stated by equation~\eqref{eq:noisyProb}. The problem is to identify the unknown control $\vec{u}$ 
given possibly perturbed measurements of the state $\vec{y}.$ Hence, the observation of the system dynamic $\mathcal{G}(\vec{u})$ 
is perturbed
by  noise $\vecsym{\eta}\in\R^K$. The noise  is assumed independent on the control $\vec{u}\in\R^d$ and normally distributed with zero mean and known covariance matrix $\vecsym{\Gamma}^{-1} \in \R^{K\times K}$, i.e. $\vecsym{\eta}\sim\mathcal{N}(0,\vecsym{\Gamma}^{-1})$. We consider a number  $J$ of ensembles (realizations of the control) combined in $\vec{U}=\left\{\vec{u}^{j} \right\}_{j=1}^J$. The EnKF is originally posed as a discrete iteration on $\vec{U}.$ The iteration index is denoted by $n$ 
and the  collection of the ensembles by $\vec{u}^{j,n}\in\R^d$, $\forall\,j=1,\dots,J$ and $n\geq 0$. According to~\cite{schillingsstuart2017}, the EnKF iterates each component of $\vec{U}^n$ at iteration $n+1$ as 
\begin{equation} \label{eq:updateEnKF}
\begin{aligned}
\vec{u}^{j,n+1} &= \vec{u}^{j,n} + \vec{C}(\vec{U}^n) \left( \vec{D}(\vec{U}^n) + \frac1{\Delta t} \vecsym{\Gamma}^{-1} \right)^{-1} (\vec{y}^{{j,n+1}} - \G(\vec{u}^{j,n}) ) \\
{\vec{y}^{j,n+1}} &= \vec{y} + \vecsym{\xi}^{j,n+1}
\end{aligned}
\end{equation}
for each $j=1,\dots,J$. Here, each observation or measurement ${\vec{y}^{j,n+1}}\in\R^K$ has been 
{perturbed by $\vecsym{\xi}^{j,n+1}\sim \mathcal{N}(0,\Delta t^{-1}\vecsym{\Sigma})$}, and $\Delta t\in\R^+$ is a parameter. As in~\cite{schillingsstuart2017} two cases for the covariance $\vecsym{\Sigma}$ will be discussed: $\vecsym{\Sigma}=0$ corresponding to a problem where the measurement data
$\vec{y}$ is unperturbed and $\vecsym{\Sigma}=\vecsym{\Gamma}^{-1}$ corresponding to the 
case where $\vecsym{\xi}^{j,n+1}$ are realizations of the noise $\vecsym{\eta}.$

Note that the update~\eqref{eq:updateEnKF} of the ensembles requires the knowledge of the operators $\vec{C}(\vec{U}^n)$ and $\vec{D}(\vec{U}^n)$ which are the covariance matrices depending on the ensemble set $\vec{U}^n$ at iteration $n$ and on $\G(\vec{U}^n)$, i.e. the image of $\vec{U}^n$ at iteration $n$. More precisely,
\begin{equation} \label{eq:covariance}
\begin{aligned}
\vec{C}(\vec{U}^n) &= \frac{1}{J} \sum_{k=1}^J \left(\vec{u}^{k,n}-\overline{\vec{u}}^n\right) \otimes \left(\G(\vec{u}^{k,n})-\overline{\G}^n\right) \in \R^{d\times K} \\
\vec{D}(\vec{U}^n) &= \frac{1}{J} \sum_{k=1}^J \left(\G(\vec{u}^{k,n})-\overline{\G}^n\right) \otimes \left(\G(\vec{u}^{k,n})-\overline{\G}^n\right) \in \R^{K\times K}
\end{aligned}
\end{equation}
where we define by $\overline{\vec{u}}^n$ and $\overline{\G}^n$ the mean of $\vec{U}^n$ and $\G(\vec{U}^n)$, namely
$$
\overline{\vec{u}}^n = \frac{1}{J} \sum_{j=1}^J \vec{u}^{j,n}, \quad \overline{\G}^n = \frac{1}{J} \sum_{j=1}^J \G(\vec{u}^{j,n}).
$$


In recent years, the EnKF was also studied as technique to solve classical and Bayesian inverse problems. For instance see the works~\cite{iglesiaslawstuart2013} and~\cite{ernstetal2015}, respectively, and the references therein. Here, we keep the attention on this type of application. The analysis of the method is proved to have a comparable accuracy with  traditional least-squares approaches to inverse problems~\cite{iglesiaslawstuart2013}. Moreover, it is known that the method provides an estimate of the unknown control $\vec{u}$ which lies in the subspace spanned by the initial ensemble set $\vec{U}^0$~\cite{iglesiaslawstuart2013}. We will see in this section that this property is still true at the continuous time level~\cite{schillingsstuart2017}. Concerning Bayesian inverse problems, instead, the method is proved to approximate specific Bayes linear estimators but it is able to provide only an approximation of the posterior measure by a (possibly weighted) sum of Dirac masses. For a detailed discussion we refer to~\cite{Apteetal2007,ernstetal2015,leglandmonbettran2009}.

As showed in~\cite{schillingsstuart2017}, it is straightforward to compute the continuous time limit equation of the update~\eqref{eq:updateEnKF}  in the general case of a nonlinear model $\G$, even if the asymptotic analysis was performed in the easier linear setting. Consider the parameter $\Delta t$ as an artificial time step for the iteration in~\eqref{eq:updateEnKF}, i.e. we take $\Delta t \sim N_t^{-1}$ where $N_t$ is the maximum number of iterations. Assume then $\vec{U}^n \approx \vec{U}(n\Delta t)=\left\{\vec{u}^{j}(n\Delta t) \right\}_{j=1}^J$ for $n\geq 0$. Scaling by $\Delta t$ and computing the limit $\Delta t\to 0^+$, the continuous time limit equation of~\eqref{eq:updateEnKF} reads 
\begin{equation} \label{eq:continuousEnKF1}
{\mathrm{d}} \vec{u}^j = \vec{C}(\vec{U}) \vecsym{\Gamma} \left( \vec{y} - \G(\vec{u}^j) \right) \, {\mathrm{dt}}  {+ \vec{C}(\vec{U}) \vecsym{\Gamma} \sqrt{\vecsym{\Sigma}} \; {\mathrm{d}\vec{W}^j}}
\end{equation}
for $j=1,\dots,J$,  initial condition $\vec{U}(0) = \vec{U}^0$ {and  ${\mathrm{d}}\vec{W}^j$ are Brownian motions}. Using the definition of the operator $\vec{C}(\vec{U})$, see~\eqref{eq:covariance}, system~\eqref{eq:continuousEnKF1} can be restated as
\begin{equation} \label{eq:continuousEnKF2}
{\mathrm{d}}  \vec{u}^j = \frac{1}{J} \sum_{k=1}^J \left\langle \G(\vec{u}^k) - \overline{\G} , \vec{y} - \G(\vec{u}^j) \right\rangle_{\vecsym{\Gamma}^{-1}} (\vec{u}^k - \overline{\vec{u}}) \, {\mathrm{dt}}  {+ \vec{C}(\vec{U}) \vecsym{\Gamma} \sqrt{\vecsym{\Sigma}} \; {\mathrm{d\vec{W}^j}}} 
\end{equation}
for $j=1,\dots,J$, where $\langle \cdot,\cdot \rangle_{\vecsym{\Gamma}^{-1}} = \langle \vecsym{\Gamma}^{\frac12} \cdot,\vecsym{\Gamma}^{\frac12} \cdot \rangle$ and $\langle \cdot,\cdot \rangle$ is the inner-product on $\R^K$.  From~\eqref{eq:continuousEnKF2} it is easy to observe that the invariant subspace property holds also at the continuous time level in 
the case $\vecsym{\Sigma} \equiv 0$ since the vector field is in the linear span of the ensemble itself.

In~\cite{schillingsstuart2017} the asymptotic behavior of the continuous time equation is analyzed in the linear setting with $\vecsym{\Sigma}\equiv 0$ so that~\eqref{eq:continuousEnKF1} is written as gradient descent equation. In fact, let us consider the case of $\G$ linear, i.e. $\G(\vec{u})=G \vec{u}$. Then the computation of the operator $\vec{C}(\vec{U})$ is
$
\vec{C}(\vec{U}) = \frac{1}J \sum_{k=1}^J \left(\vec{u}^k-\overline{\vec{u}}\right) \left(\vec{u}^k-\overline{\vec{u}}\right)^T G^T.
$
Further, note that the least squares functional~\eqref{eq:leastSqFnc} yields
\begin{equation} \label{eq:gradientGLinear}
\nabla_\vec{u} \Phi(\vec{u},\vec{y}) = - G^T \vecsym{\Gamma} ( \vec{y} - G \vec{u} ).
\end{equation}
Therefore, equation~\eqref{eq:continuousEnKF1} is stated in terms of the gradient of $\Phi$ as 
\begin{equation} \label{eq:gradientEq}
\frac{\mathrm{d}}{\mathrm{d}t} \vec{u}^j = - \frac{1}J \sum_{k=1}^J (\vec{u}^k-\overline{\vec{u}}) \otimes ( \vec{u}^k - \overline{\vec{u}} )  \nabla_\vec{u} \Phi(\vec{u}^j,\vec{y})
\end{equation}
for $j=1,\dots,J$. Equation~\eqref{eq:gradientEq} describes a preconditioned gradient descent equation for each ensemble. In fact, $\vec{C}(\vec{U})$ is positive semi-definite and hence
$$
\frac{\mathrm{d}}{\mathrm{d}t} \Phi(\vec{u}(t),\vec{y}) = \frac{\mathrm{d}}{\mathrm{d}t} \frac12 \left\| \vecsym{\Gamma}^{\frac12} \left(\vec{y}-G\vec{u}\right) \right\|^2 \leq 0.
$$
Observe that, although the forward operator is assumed to be linear, the gradient flow is nonlinear. For further details and properties of the gradient descent equation~\eqref{eq:gradientEq} we refer to~\cite{schillingsstuart2017}. In particular, here we recall the important result on the velocity of the collapse of the ensembles towards their mean in the large time limit.
\begin{lemma}[Theorem 3 in~\cite{schillingsstuart2017}] \label{lem:schillingsstuart}
	Let $\vec{U}^0$ be the initial set of ensembles. Then the matrix $\vec{R}(t)$ whose entries are
	$$
	\left( \vec{R}(t) \right)_{ij} = \left\langle G(\vec{u}^i-\overline{\vec{u}}),G(\vec{u}^j-\overline{\vec{u}}) \right\rangle_{\vecsym{\Gamma}}
	$$
	converges to $0$ for $t\to\infty$ and indeed $\left\| \vec{R}(t) \right\| = O(Jt^{-1})$.
\end{lemma}
The previous Lemma also states that the collapse slows down linearly as the ensemble size increases. Later, this property is also obtained in the mean-field limit for a large ensemble size.

We point out that the continuous time limit derivation suggests to stop at time $t = 1$. However, the study of the long-time analysis of the ODE system, such as the stability analysis in Section~\ref{sec:unstableMoments} and Section~\ref{sec:stableAnalysis}, highlights possible improvements of the
algorithm.

\section{Mean-field limit of the Ensemble Kalman Filter} \label{sec:meanfield}

Typically, the EnKF method is applied for a fixed and finite ensemble size. In fact, it is clear from~\eqref{eq:updateEnKF} and~\eqref{eq:continuousEnKF2} that the computational and memory cost of the method increases with the number of the ensembles. The analysis of the method was also studied in the large ensemble limit, see e.g.~\cite{ernstetal2015,kwiatkowskimandel2015,lawtembinetempone2016,leglandmonbettran2009}.
However, to the best of our knowledge, the derivation of a kinetic equation that holds in the limit of a large 
number of ensembles has not yet been proposed. In this section, we derive the corresponding mean-field limit of the continuous time equation focusing on the case of a linear model $G$ and with $\vecsym{\Sigma}=\vec{0}$ as in \cite{schillingsstuart2017}. 

We follow the classical formal derivation to formulate a mean-field equation of a particle system, see~\cite{CarrilloFornasierToscaniVecil2010,hatadmor2008,PareschiToscaniBOOK,Toscani2006}. Let us denote by
\begin{equation} \label{eq:kineticf}
f = f(t,\vec{u}) : \R^+ \times \R^d \to \R^+
\end{equation}
the compactly supported on $\R^d$ probability density of $\vec{u}$ at time $t$ and introduce the first moment $\vec{m}\in\R^d$ and the second moment $\vec{E}\in\R^{d\times d}$ of $f$ at time $t$, respectively, as
\begin{equation} \label{eq:moments}
\vec{m}(t) = \int_{\R^d} \vec{u} f(t,\vec{u}) \mathrm{d}\vec{u}, \quad \vec{E}(t) = \int_{\R^d} \vec{u} \otimes \vec{u} f(t,\vec{u}) \mathrm{d}\vec{u}.
\end{equation}
Since $\vec{u}\in\R^d$, the corresponding discrete measure on the ensemble set $\vec{U} = \left\{ \vec{u}^j \right\}_{j=1}^J$ is therefore given  by the empirical measure 
\begin{equation} \label{eq:empiricalf}
f(t,\vec{u}) = \frac{1}J \sum_{j=1}^J \delta(\vec{u}^j - \vec{u}) = \frac{1}J \sum_{j=1}^J \prod_{i=1}^d \delta(u^j_i - u_i),
\end{equation}
where $u^j_i\in\R$ is the component $i$ of the $j$-th ensemble. Let us define the operator
$$
	\vecsym{\C}(\vec{U}) = \frac{1}J \sum_{k=1}^J (\vec{u}^k-\overline{\vec{u}}) \otimes (\vec{u}^k-\overline{\vec{u}})
$$
with the corresponding entry 
$$
	\left( \vecsym{\C}(\vec{U}) \right)_{\kappa,\ell} = \frac{1}J \sum_{k=1}^J u_\kappa^k u_\ell^k - \overline{u}_\kappa \frac{1}J \sum_{k=1}^J u_\ell^k - \overline{u}_\ell \frac{1}J \sum_{k=1}^J u_\kappa^k + \overline{u}_\kappa \overline{u}_\ell = \frac{1}J \sum_{k=1}^J u_\kappa^k u_\ell^k - \overline{u}_\kappa \overline{u}_\ell,
$$
where $\overline{u}_i$ denotes the component $i$ of the mean $\overline{\vec{u}}$ of the ensembles. This formulation
allows for a mean-field limit as 
$$
\left(\vecsym{\C}(t)\right)_{\kappa,\ell} = \int_{\R^d} u_\kappa u_\ell f(t,\vec{u}) \mathrm{d}\vec{u} - \int_{\R^d} u_\kappa f(t,\vec{u}) \mathrm{d}\vec{u} \int_{\R^d} u_\ell f(t,\vec{u}) \mathrm{d}\vec{u}
$$
and therefore $\vecsym{\C}(\vec{U})$ can be written in terms of the moments~\eqref{eq:moments} of the empirical measure only as
\begin{equation} \label{eq:covarianceMeanField}
\vecsym{\C}(t) = \vec{E}(t) - \vec{m}(t) \otimes \vec{m}(t).
\end{equation}
Let us denote $\varphi(\vec{u}) \in C_0^1(\R^d)$ a sufficiently smooth test function. We compute
\begin{align*}
\frac{\mathrm{d}}{\mathrm{d}t} \left\langle f , \varphi \right\rangle &= \frac{\mathrm{d}}{\mathrm{d}t} \int_{\R^d} \frac{1}{J} \sum_{j=1}^J \delta(\vec{u} - \vec{u}^j) \varphi(\vec{u}) \mathrm{d}\vec{u} = - \frac{1}{J} \sum_{j=1}^J \nabla_\vec{u} \varphi(\vec{u}^j) \cdot \vecsym{\C}(t) \nabla_\vec{u} \Phi(\vec{u}^j,\vec{y}) \\
&= - \int_{\R^d} \nabla_\vec{u} \varphi(\vec{u}) \cdot \vecsym{\C}(t) \nabla_\vec{u} \Phi(\vec{u},\vec{y}) f(t,\vec{u}) \mathrm{d}\vec{u}
\end{align*}
which finally leads to the following strong form of the mean-field kinetic equation corresponding to the gradient descent equation~\eqref{eq:gradientEq}:
\begin{equation} \label{eq:kineticFromEnKF}
\partial_t f(t,\vec{u}) - \nabla_\vec{u} \cdot \left( \vecsym{\C}(t) \nabla_\vec{u} \Phi(\vec{u},\vec{y}) f(t,\vec{u}) \right) = 0.
\end{equation}

Equation~\eqref{eq:kineticFromEnKF} provides a closed formula for the evolution in time of the distribution $f$ of the unknown control $\vec{u}$ when the observations $\vec{y}$ and the linear model $G$ are given and when endowed with an initial guess $f^0(\vec{u}) = f(t=0,\vec{u})$ for the unknown control.


\subsection{Moment equations and linear stability analysis}
\label{sec:unstableMoments}

As discussed in Section~\ref{sec:enkf}, the EnKF computes a solution to the inverse problem as mean of the ensembles in the large time behavior. Since the kinetic equation~\eqref{eq:kineticFromEnKF} formally holds in the limit of a large number of ensembles, here we analyze approximations to the solution of the inverse problem provided by the first moment $\vec{m}(t)$ of the kinetic distribution, see~\eqref{eq:moments}.

Due to definition~\eqref{eq:moments}, multiplying~\eqref{eq:kineticFromEnKF} by $\vec{u}$, integrating over $\R^d$ and integrating by parts the second term, we get the following evolution equation for the first moment:
$$
\frac{\mathrm{d}}{\mathrm{d}t} \vec{m}(t) + \int_{\R^d} \vecsym{\C}(t) \nabla_\vec{u} \Phi(\vec{u},\vec{y}) f(t,\vec{u}) \mathrm{d}\vec{u} = \vec{0}.
$$
In particular, since we are assuming the simple setting of a linear model $\G(\vec{u}) = G\vec{u}$, using~\eqref{eq:gradientGLinear}, we can explicitly compute the integral and obtain
\begin{equation} \label{eq:linear1stMomEq}
\frac{\mathrm{d}}{\mathrm{d}t} \vec{m}(t) + \vecsym{\C}(t)  \nabla_\vec{u} \Phi(\vec{\vec{m}},\vec{y}) = \vec{0}.
\end{equation}
Multiplying~\eqref{eq:kineticFromEnKF} by $\vec{u} \otimes \vec{u}$ and integrating over $\R^d$  we obtain the following evolution equation for the second moment:
\begin{equation} \label{eq:linear2ndMomEq}
\frac{\mathrm{d}}{\mathrm{d}t} \vec{E}(t) + \sum_{k=1}^d \int_{\R^d} \vec{T}_k^{(1)}(\vec{u}) \left( \vecsym{\C}(t) \nabla_\vec{u} \Phi(\vec{u},\vec{y}) f(t,\vec{u}) \right)_k \mathrm{d}\vec{u} = \vec{0},\quad \vec{T}_k^{(1)}(\vec{u}) = \frac{\partial}{\partial u_k} \vec{u} \otimes \vec{u}.
\end{equation}

Hence, equation~\eqref{eq:linear1stMomEq} and equation~\eqref{eq:linear2ndMomEq} provide a closed system of ordinary differential equations.

\begin{remark} \label{rem:mSolution}
	As in Bayesian approach to inverse problems, also equation~\eqref{eq:kineticFromEnKF} poses the problem of selecting a solution out of $f$ which only provides a distribution for the unknown control $\vec{u}$. As pointed out at the beginning of this subsection, since the kinetic equation is derived via mean-field limit we choose, accordingly to the solution provided by the EnKF, the expected value $\vec{m}$ as an estimator of the unknown parameter $\vec{u}$. Observe that a steady-state $\vec{m}^\infty$ of equation~\eqref{eq:linear1stMomEq} is given by
	$$
	\vec{m}^\infty = \arg\min_{\vec{u}} \Phi(\vec{u},\vec{y}),
	$$
	corresponding to a control that minimizes the least squares functional $\Phi$.  In the case of a linear model $G$, the above condition can be also stated as $\vec{y}-G\vec{u} \in \ker G^T.$ Neither $\vec{u}$ nor $\vec{m}^\infty$ need to be unique. 
\end{remark}

Equation~\eqref{eq:linear1stMomEq} for the first moment $\vec{m}$ and~\eqref{eq:linear2ndMomEq} for the second moment $\vec{E}$ give rise to a coupled system of ordinary differential equations. In the following, we employ a stability analysis for these equations in the simple case of a one-dimensional control in order to analyze the stability of the estimator $\vec{m}$. 

First, we observe that in the case of a scalar control the system of the moment equations reduces to
\begin{equation} \label{eq:system1D}
\begin{aligned}
\frac{\mathrm{d}}{\mathrm{d}t} m(t) &= G ( E(t) - m^2(t) )(y - G m(t)) \\
\frac{\mathrm{d}}{\mathrm{d}t} E(t) &= 2 G ( E(t) - m^2(t) ) (y m(t) - G E(t))	
\end{aligned}
\end{equation}
with $y\in\R$ and $G\in \R\setminus\{0\}$. The nullclines of the system of ODEs~\eqref{eq:system1D} are given by
$$
m = \frac{y}{G}, \quad E = \frac{y}{G} m, \quad E = m^2.
$$
The equilibrium or fixed points of~\eqref{eq:system1D} are the intersections of the nullclines and therefore we have the following three sets of points:
$$
F_0 = (0,0), \quad F_1 = (\frac{y}{G},\frac{y^2}{G^2}), \quad F_k=(k,k^2), \; k\in\R,
$$
i.e. all the fixed points are on the parabola $E = m^2$ in the phase plane $(m,E)$. Given the Jacobian $\vec{J}\in\R^{2\times 2}$ of the ODE system~\eqref{eq:system1D}
\begin{equation} \label{eq:jacobianODE}
\vec{J}(m,E) = \begin{bmatrix}
3 G^2 m^2 - 2 G y m - G^2 E & -G^2 m + G y \\
&\\
2 G y E + 4 G^2 m E - 6 G y m^2 & -4 G^2 E + 2 G y m + 2 G^2 m^2 
\end{bmatrix}	
\end{equation}
it follows that $\vec{J}(F_k)$ has eigenvalues $\mu_1 = \mu_2 = 0$. Clearly, the same holds for $F_0$ and $F_1$, since they are points of the type $F_k$. Therefore all the fixed points are non-hyperbolic and the stability must be analyzed directly. More precisely, since $\mu_1 = \mu_2 = 0$, the fixed points are Bogdanov-Takens-type equilibria and hence unstable as we indeed show in the following analysis.

\begin{figure}[t!]
	\centering
	\includegraphics[width=0.49\textwidth]{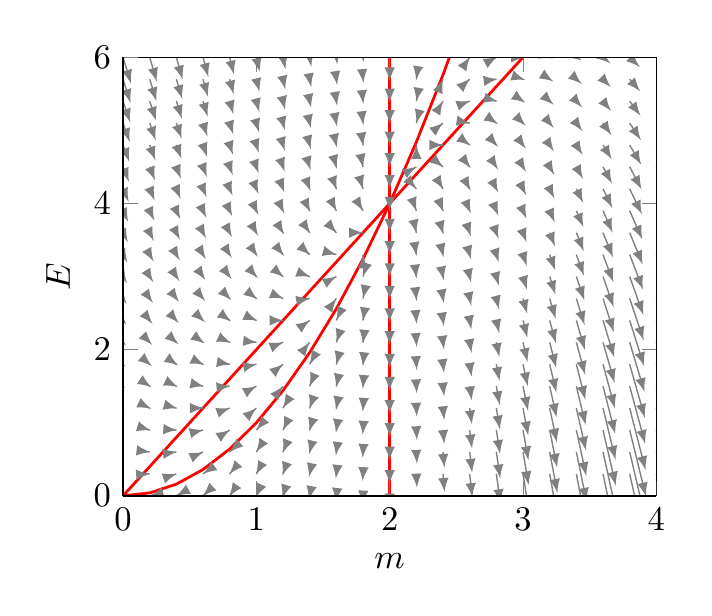}
	\includegraphics[width=0.49\textwidth]{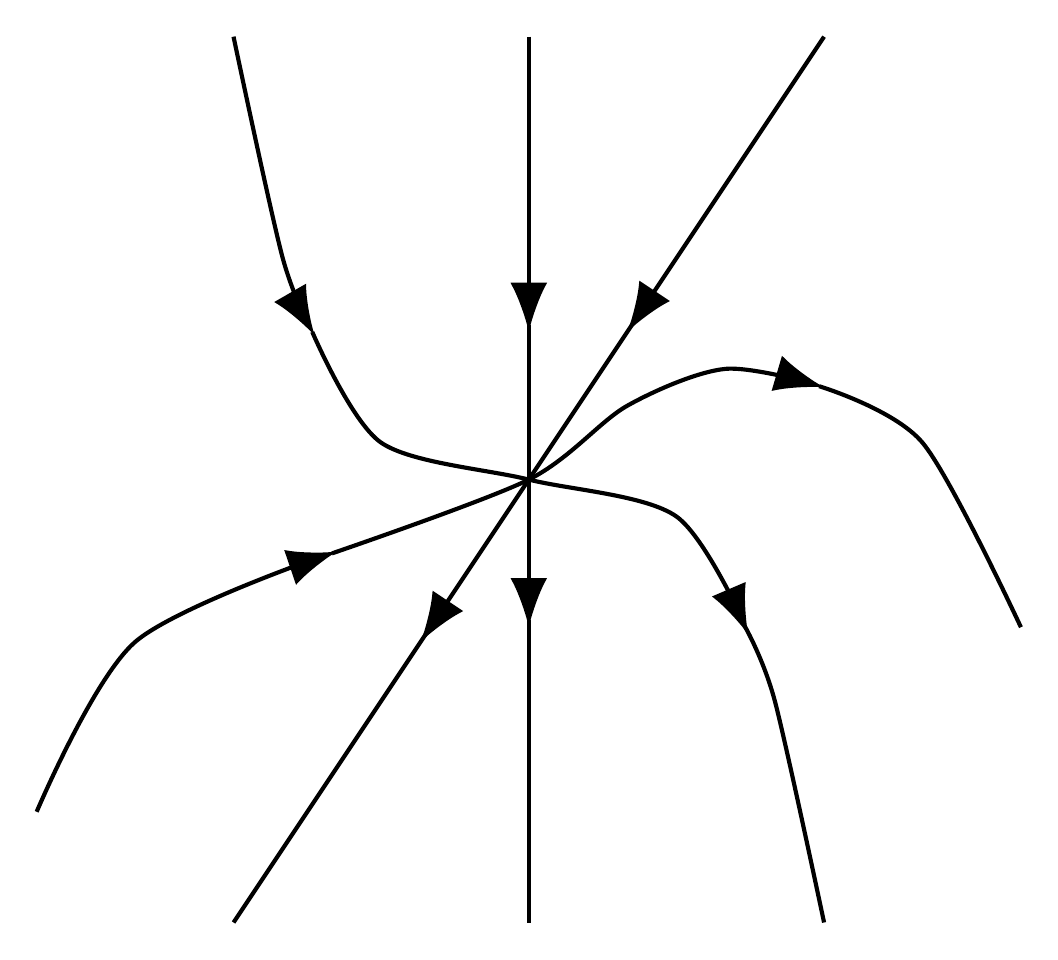}
	\caption{Left: vector field of the ODE system~\eqref{eq:system1D} with $(y,G)=(2,1)$. Red lines are the nullclines. Right: trajectory behavior around the equilibrium $(\frac{y}{G},\frac{y^2}{G^2}) = (2,4)$.\label{fig:phaseplane}}
\end{figure}

The vector field of the system~\eqref{eq:system1D} can be easily analyzed on nullclines and on the $m$- and $E$-axis of the phase plane. For the sake of simplicity, let us assume that $\frac{y}{G} > 0$. The analysis is equivalent in the opposite case. Let $m(0) = \frac{y}{G}$ so that $\frac{\mathrm{d}}{\mathrm{d}t} m=0$ for all $t$. We have that
$$
\frac{\mathrm{d}}{\mathrm{d}t} E = - \frac{2}{G^2} (y^2 - G^2 E )^2 < 0
$$
and therefore $E$ is decreasing in time on the nullcline $m = \frac{y}{G}$ which in turn means that $F_1$ is an attractor only if $E(0) > \frac{y}{G^2}$. Let now $E(0) = \frac{y}{G} m(0)$, for some $m(0)$. Then we have $\frac{\mathrm{d}}{\mathrm{d}t} E = 0$ and
$$
\frac{\mathrm{d}}{\mathrm{d}t} m = m (y - G m)^2 = \begin{cases}
> 0, & \text{if $m(0)>0$},\\
< 0, & \text{otherwise}.
\end{cases}
$$ 
Thus, since $m(0)>0$ is the only acceptable initial condition in order to guarantee that $E(0)>0$, the trajectories are moving on the right side of the phase plane on the nullcline $E = \frac{y}{G} m$. Obviously, each trajectory is still in time on the nullcline $E = m^2$ since $\frac{\mathrm{d}}{\mathrm{d}t} m = \frac{\mathrm{d}}{\mathrm{d}t} E = 0$. The nullclines and the complete vector field for the case $(y,G)=(2,1)$ is shown in the left panel of Figure~\ref{fig:phaseplane}. We immediately observe that the behavior around the equilibrium point $F_1$ is unstable as showed also in the right panel of Figure~\ref{fig:phaseplane}.

The previous considerations can be also derived by looking at the solutions of~\eqref{eq:system1D}. Assuming that the initial conditions are such that $E(0) \neq m(0)^2$, we get then the following pairs of analytical solutions:
\begin{gather*}
m(t) = \frac{y}{G}, \quad E(t) = \frac{y^2}{G^2} + \frac{1}{2G^2(C+t)}\\
m(t) = \frac{y}{G} \pm \frac{1}{G\sqrt{-2 C_1 G t - 2 C_2 G}}, \quad E(t) = m^2 + \frac{\frac{\mathrm{d}}{\mathrm{d}t} m}{G(y - G m)}
\end{gather*}
with $C, C_1, C_2 \in \R$ constants uniquely prescribed by the initial conditions. In particular, the first set of solutions is found by assuming that $m(0) = \frac{y}{G}$ and solving the following Riccati's equation with constant coefficients
$$
\frac{\mathrm{d}}{\mathrm{d}t} E(t) = -2 G^2 E^2 + 4 E y - 2 \frac{y^4}{G^2}.
$$
In this case, letting $E(0)=E_0$, the constant is given by $C = \frac12 (G^2 E_0^2 - y^2)$ which is positive when $E_0 > \frac{y^2}{G^2}$ and negative otherwise. In this latter case we also observe that there exists a time $t$ in which the trajectory $E(t)$ has a vertical asymptote. For the above discussion we know that $E(t)$ is also decreasing. It is also simple to observe that in the second pair of solutions $E(t)$ can blow up driving $m(t)$ away from the equilibrium.

\begin{remark}
	We observe that the linear stability analysis provided in this section is applied to system~\eqref{eq:system1D} without including restrictions on $E(t)-m^2(t)$ which must be non-negative. Under this constraint, the region $E(t)<m^2(t)$ is not admissible and the unstable equilibrium $(\frac{y}{G},\frac{y^2}{G^2})$ lies on the boundary of this region.
\end{remark}


\section{Extension of the mean-field EnKF method} \label{sec:stabilization}

The analysis of Section~\ref{sec:unstableMoments} shows that, at least in a one-dimensional setting, the system of moment equations~\eqref{eq:system1D} could lead to unconditionally unstable equilibria. This is due to the possible decay of the energy which drives the expected value far from the equilibrium value. In the general case of a $d$-dimensional control $\vec{u}$, the situation may be even more complex.

The instability of  fixed points of~\eqref{eq:system1D} can be related to the loss of an $O(\Delta t)$ term in the derivation of the continuous time limit equation~\eqref{eq:continuousEnKF1}. In fact, instability could occur also for~\eqref{eq:gradientEq} but it is possible to show that the discrete equation~\eqref{eq:updateEnKF} has stable equilibria.

Next,  we stabilize the system of the moment equations~\eqref{eq:system1D} by introducing additional uncertainty to the microscopic interactions. This leads to a diffusive term in the kinetic equation avoiding the decay  of  kinetic energy and the appearance of unstable equilibria. First, we write binary microscopic interactions corresponding to the mean-field kinetic equation~\eqref{eq:kineticFromEnKF}. Then, we introduce noise  in these interactions and we derive a Fokker-Planck-type equation. Finally, we study the stability of the resulting moment system. 

Let again $f=f(t,\vec{u}):\R^+ \times \R^d \to \R$ be the probability density of the control $\vec{u}\in\R^d$ at time $t>0$ as defined in~\eqref{eq:kineticf}. Let $\vec{m}\in\R^d$ and $\vec{E}\in\R^{d\times d}$ be the first and the second moment of $f$, respectively, as given in~\eqref{eq:moments}. We introduce the  microscopic interaction rules:
\begin{equation} \label{eq:microInteraction}
\begin{aligned}
\vec{u} &= \vec{u}_* - \epsilon (\vec{E} - \vec{m} \otimes \vec{m}) \nabla_\vec{u} \Phi(\vec{u}_*,\vec{y}) + \sqrt{\epsilon} \, \vec{K}(\vec{u}_*) \vecsym{\xi} \\ &= \vec{u}_* - \epsilon \, \vecsym{\C}(t) \nabla_\vec{u} \Phi(\vec{u}_*,\vec{y}) + \sqrt{\epsilon} \, \vec{K}(\vec{u}_*) \vecsym{\xi}
\end{aligned}
\end{equation}
where $\vec{u}$ is the post-interaction value of the ensemble member, $\vec{u}_*$ is its pre-interaction value and $\vecsym{\xi}\in\R^d$ is a random variable with given distribution $\theta(\vecsym{\xi})$ having zero mean and covariance matrix $\vecsym{\Lambda}\in\R^{d\times d}$. Instead, $\vec{K}(\vec{u}_*)\in\R^{d\times d}$ is an arbitrary function of $\vec{u}_*$. For $\vec{K} = \vec{C}$ we observe a similar structure as in equation~\eqref{eq:continuousEnKF1}. The quantity $\epsilon$ describes the strength of the interactions and it is a scattering rate.

\begin{remark} \label{rem:probabilisticMicro}
	Observe that~\eqref{eq:microInteraction} is in fact the microscopic interaction corresponding to the mean-field equation~\eqref{eq:kineticFromEnKF}, that is in the case of a linear model
	\begin{equation} \label{eq:probabInteraction}
	\vec{u} = \left( \vec{1} - \epsilon \, \vecsym{\C}(t) G^T \vecsym{\Gamma} G \right) \vec{u}_* + \epsilon \, \vecsym{\C}(t) G^T \vecsym{\Gamma} G G^{-1} \vec{y},
	\end{equation}
	with an additional term representing the uncertainty in the interaction. The interaction~\eqref{eq:probabInteraction} has a probabilistic interpretation~\cite{AlbiPareschi2013} provided
	$$
	\epsilon \, \rho\left(\vecsym{\C}(t) G^T \vecsym{\Gamma} G \right) \leq 1,
	$$
	where $\rho(\cdot)$ is the spectral radius.
\end{remark}
The probability density $f$ satisfies the following (linear) Boltzmann equation in weak form
\begin{equation} \label{eq:wfBoltzmann}
\frac{\mathrm{d}}{\mathrm{d}t} \int_{\R^d} f(t,\vec{u}) \varphi(\vec{u}) \mathrm{d}\vec{u} = \left\langle \int_{\R^d} (\varphi(\vec{u})-\varphi(\vec{u}_*)) f(t,\vec{u}) \mathrm{d}\vec{u} \right\rangle
\end{equation}
where $\varphi \in C_c^\infty(\R^d)$ is a test function and where the operator $\langle \cdot \rangle$ denotes the mean with respect to the distribution $\theta$, i.e. $
\langle g  \rangle = \int_{\R^d} g(\vecsym{\xi}) \theta(\vecsym{\xi}) \mathrm{d}\vecsym{\xi}.
$
Consider the time asymptotic scaling by setting
\begin{equation} \label{eq:scaling}
\tau = t \epsilon, \quad f(t,\vec{u}) = \tilde{f}(\tau,\vec{u})
\end{equation}
and allow $\epsilon \to 0^+$. This corresponds to large interaction frequencies and small interaction strengths, a situation similar to the so-called grazing collision limit ~\cite{Desvillettes,DiPernaLions,PareschiToscaniVillani,Villani1999}. We denote the scaled quantities again by $f$ and $t,$ respectively. A second-order Taylor expansion yields the corresponding formal Fokker-Planck equation: 
\begin{align*}
\varphi(\vec{u}) - \varphi(\vec{u}_*) 
=& \nabla_\vec{u} \varphi(\vec{u}_*) \cdot (\vec{u}-\vec{u}_*) + \frac12 (\vec{u}-\vec{u}_*)^T \vec{H}(\varphi(\vec{u}_*)) (\vec{u}-\vec{u}_*) \\ &+ \frac12 (\vec{u}-\vec{u}_*)^T \widehat{\vec{H}}(\varphi;\tilde{\vec{u}},\vec{u}_*) (\vec{u}-\vec{u}_*)
\end{align*}
with $\tilde{\vec{u}} = \alpha \vec{u}_* + (1-\alpha) \vec{u}$, $\alpha \in (0,1)$ and where $\vec{H}\in\R^{d\times d}$ is the Hessian matrix and
$
\widehat{\vec{H}}(\varphi;\tilde{\vec{u}},\vec{u}_*) = \vec{H}(\varphi(\tilde{\vec{u}})) - \vec{H}(\varphi(\vec{u}_*)).
$
Substituting this expression in equation~\eqref{eq:wfBoltzmann} and using definition~\eqref{eq:microInteraction} of the microscopic interactions, we obtain 
\begin{align*}
& \frac{\mathrm{d}}{\mathrm{d}t} \int_{\R^d} f(t,\vec{u}) \varphi(\vec{u}) \mathrm{d}\vec{u} = \frac{1}{\epsilon} \left( \mathcal{A} + \mathcal{B} + \mathcal{R} \right),     \\
& \mathcal{A} = - \epsilon \left\langle \int_{\R^d} \nabla_\vec{u} \varphi(\vec{u}_*) \cdot \left( \vecsym{\C}(t) \nabla_\vec{u} \Phi(\vec{u}_*,\vec{y}) \right) f(t,\vec{u}_*) \mathrm{d}\vec{u}_* \right\rangle 
%
+ \sqrt{\epsilon} \left\langle \int_{\R^d} \nabla_\vec{u} \varphi(\vec{u}_*) \cdot \vec{K}(\vec{u}_*) \vecsym{\xi} f(t,\vec{u}_*) \mathrm{d}\vec{u}_* \right\rangle, \\
& \mathcal{B} = \frac{\epsilon^2}{2} \left\langle \int_{\R^d} \left( \vecsym{\C}(t) \nabla_\vec{u} \Phi(\vec{u}_*,\vec{y}) \right)^T \vec{H}(\varphi(\vec{u}_*)) \left( \vecsym{\C}(t) \nabla_\vec{u} \Phi(\vec{u}_*,\vec{y}) \right) f(t,\vec{u}_*) \mathrm{d}\vec{u}_* \right\rangle \\
&- \epsilon \sqrt{\epsilon} \left\langle \int_{\R^d} \left(\vec{K}(\vec{u}_*)\vecsym{\xi}\right)^T \vec{H}(\varphi(\vec{u}_*)) \left( \vecsym{\C}(t) \nabla_\vec{u} \Phi(\vec{u}_*,\vec{y}) \right) f(t,\vec{u}_*) \mathrm{d}\vec{u}_* \right\rangle \\
%
&+ \frac{\epsilon}{2} \left\langle \int_{\R^d} \Tr\left( (\vecsym{\xi}\otimes\vecsym{\xi})^T \vec{K}(\vec{u}_*)^T \vec{H}(\varphi(\vec{u}_*)) \vec{K}(\vec{u}_*) \right) f(t,\vec{u}_*) \mathrm{d}\vec{u}_* \right\rangle,
\end{align*}
where $\Tr(\cdot)$ is the matrix trace and $\mathcal{R}$ is the remaining term. One can easily prove that $\epsilon^{-1}\mathcal{R}$ vanishes in the asymptotic scaling~\eqref{eq:scaling}. In order to show this, it is sufficient the fact that $\varphi$ is an enough smooth function and thus each second partial derivative is Lipschitz continuous so that $\exists\,L>0$ such that
$$
	\left| \frac{\partial^2 \varphi(\tilde{\vec{u}})}{\partial u_iu_j} - \frac{\partial^2 \varphi(\vec{u}_*)}{\partial u_iu_j} \right| \leq L |\tilde{\vec{u}} - \vec{u}_*| < L |\vec{u} - \vec{u}_*| = L | \epsilon \vecsym{\C}(t) \nabla_\vec{u} \Phi(\vec{u}_*,\vec{y}) + \sqrt{\epsilon} \vec{K}(\vec{u}_*) \vecsym{{\xi}} | \xrightarrow{\epsilon\to 0^+} 0
$$
for all $i,j=1,\dots,d$.



For $\vec{K}(\vec{u}_*)$ constant or depending on moments of the kinetic distribution $f$, the grazing limit in strong form is then obtained as 
\begin{equation} \label{eq:kineticFP}
\partial_t f(t,\vec{u}) = \nabla_\vec{u} \cdot \left( \vecsym{\C}(t) \nabla_\vec{u} \Phi(\vec{u},\vec{y}) f(t,\vec{u}) \right) + \frac12 \nabla_\vec{u} \cdot \left( \vecsym{\Lambda} \vec{K}^T\vec{K} \nabla_\vec{u} f(t,\vec{u}) \right)
\end{equation}
where we used the basic fact
$$
\Tr\left( \vecsym{\Lambda} \vec{K}^T\vec{K} \vec{H}(f(t,\vec{u})) \right) = \nabla_\vec{u} \cdot \left( \vecsym{\Lambda} \vec{K}^T\vec{K} \nabla_\vec{u} f(t,\vec{u}) \right). 
$$

Some remarks are in order. As expected, the Fokker-Planck-type equation~\eqref{eq:kineticFP} is consistent with the kinetic equation~\eqref{eq:kineticFromEnKF} in the limit of vanishing covariance $\vecsym{\Lambda}$. The introduction of the uncertainty in~\eqref{eq:microInteraction} allows for a different interpretation of the data perturbation in~\eqref{eq:updateEnKF} and~\eqref{eq:continuousEnKF2} within the kinetic model.


\subsection{Moment equations and linear stability analysis} \label{sec:stableAnalysis}

In the setting of \cite{schillingsstuart2017} we have $\G(\vec{u})=G\vec{u}$ and, for $\vec{K}=\vec{I}$ identity matrix, a straightforward computation leads to the following moment equations
based on the Fokker-Planck equation \eqref{eq:kineticFP}. 
\begin{equation} \label{eq:stableMoments}
\begin{aligned}
\frac{\mathrm{d}}{\mathrm{d}t} \vec{m}(t) &= - \vecsym{\C}(t) \nabla_\vec{u} \Phi(\vec{m},\vec{y})\\
\frac{\mathrm{d}}{\mathrm{d}t} \vec{E}(t) &= - \sum_{k=1}^d \int_{\R^d} \vec{T}_k^{(1)}(\vec{u}) \left( \vecsym{\C}(t) \nabla_\vec{u} \Phi(\vec{u},\vec{y}) f(t,\vec{u}) \right)_k \mathrm{d}\vec{u} + \frac12 \sum_{i,j=1}^d \Lambda_{ij}^2 \int_{\R^d} \vec{T}_{ij}^{(2)}(\vec{u}) f(t,\vec{u}) \mathrm{d}\vec{u},
\end{aligned}
\end{equation}
where $ \vec{m}$ and $\vec{E}$ are defined as before and where we have 
$$
\vec{T}_k^{(1)}(\vec{u}) = \frac{\partial}{\partial u_k} \vec{u} \otimes \vec{u}, \quad \vec{T}_{ij}^{(2)}(\vec{u}) = \frac{\partial}{\partial u_i u_j} \vec{u} \otimes \vec{u}.
$$
Comparing  equation  \eqref{eq:stableMoments} and~\eqref{eq:linear1stMomEq}, we observe that they are equivalent. This implies that the equation for $\vec{m}$ is still providing a solution according to Remark~\ref{rem:mSolution}. 
Instead, the equation of the second moment $\vec{E}$ has  an additional term that stabilize the equilibria of~\eqref{eq:stableMoments}. 

\begin{figure}[t!]
	\centering
	\includegraphics[width=0.49\textwidth]{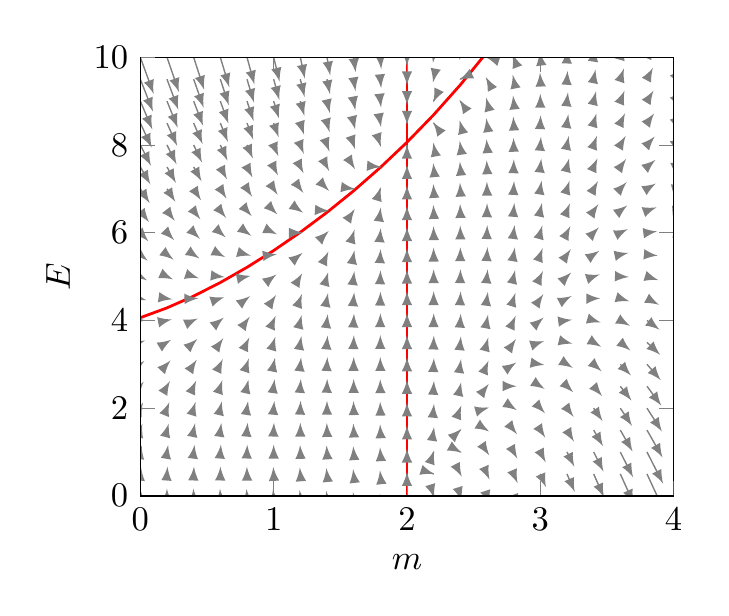}
	\includegraphics[width=0.49\textwidth]{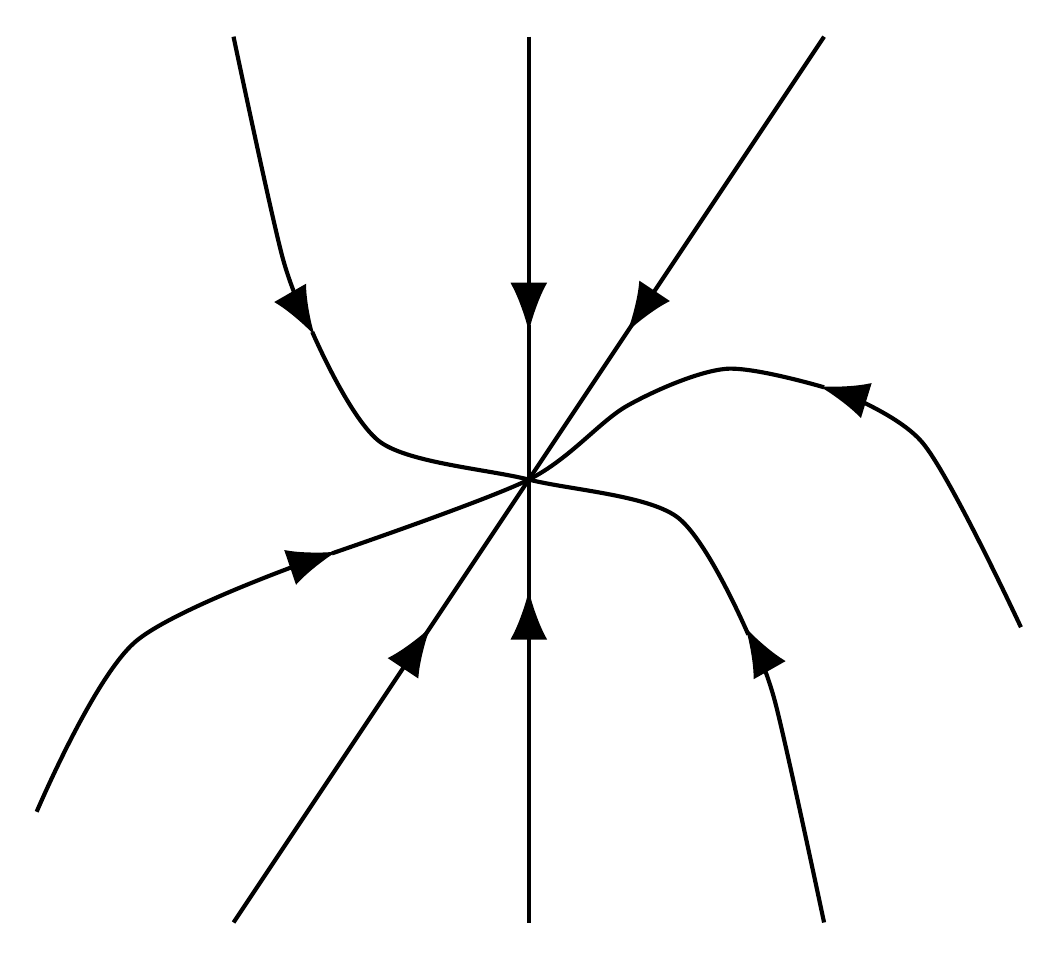}
	\caption{Left: vector field of the ODE system~\eqref{eq:systemNoise1D} with $(y,G)=(2,1)$. Red lines are the nullclines. Right: trajectory behavior around the equilibrium $\tilde{F}_1^+ = (2,8)$.\label{fig:phaseplaneNoise}}
\end{figure}

We analyze linear stability of~\eqref{eq:stableMoments} in the case of a one-dimensional control. In this particular case  the moment equations are
\begin{equation} \label{eq:systemNoise1D}
\begin{aligned}
\frac{\mathrm{d}}{\mathrm{d}t} m(t) &= G ( E(t) - m(t)^2 ) (y - G m(t)) \\
\frac{\mathrm{d}}{\mathrm{d}t} E(t) &= 2 G ( E(t) - m(t)^2 ) (y m(t) - G E(t)) + \lambda^2	
\end{aligned}
\end{equation}
with $y\in\R$, $G\in\R\setminus\{0\}$ and where now $\lambda^2\in\R$ represents the variance of the univariate noise $\xi$. Following the same analysis performed in Section~\ref{sec:unstableMoments}, we compute the nullclines of the ODE system~\eqref{eq:systemNoise1D} and they are given by
$$
m = \frac{y}{G}, \quad E = m^2, \quad E = \frac{m(y+Gm) \pm \sqrt{m^2(y-Gm)^2 + 2 \lambda^2}}{2G}.
$$
We are interested in the behavior around the equilibrium with $m=\frac{y}{G}$ which is obtained as intersection of the first and the third nullcline:
$$
\tilde{F}_1^\pm = (\frac{y}{G},\frac{y^2}{G^2} \pm \frac{\sqrt{2\lambda^2}}{2G}).
$$
Observe that this equilibrium point is in fact the equilibrium point $F_1$ given in Section~\ref{sec:unstableMoments} when $\lambda \to 0^+$. For simplicity, in the following we consider $y,G>0$. Similar considerations can be done in the other cases. Letting $m(0)=\frac{y}{G}$ so that $\frac{\mathrm{d}}{\mathrm{d}t} m = 0$ for all $t$, we have
$$
\frac{\mathrm{d}}{\mathrm{d}t} E = -2G^2 \left( E - \frac{y^2}{G^2} \right)^2 + \lambda^2
$$
where the right-hand side represents a parabola in $E$ with negative leading coefficient. Therefore, using classical arguments of stability theory for ODEs, we can state that the greater root $\tilde{F}_1^+$ is the stable equilibrium and the smaller root $\tilde{F}_1^-$ is the unstable equilibrium. This result can be also obtained by looking at the eigenvalues of the Jacobian matrix of the system~\eqref{eq:systemNoise1D} which is equivalent to~\eqref{eq:jacobianODE}. In fact, computing the eigenvalues $\mu^\pm_{1,2}$ of $\vec{J}(m,E)$ evaluated in $\tilde{F}_1^\pm$ we have
$$
\mu_1^\pm = \mp \frac{G}{2} \sqrt{2\lambda^2}, \quad \mu_2^\pm = \mp 2G\sqrt{2\lambda^2}
$$
and therefore the equilibrium $\tilde{F}_1^+$ corresponding to the two negative eigenvalues is stable. Moreover, we stress the fact that the in the case of~\eqref{eq:systemNoise1D} the equilibria are no longer non-hyperbolic as in the case of~\eqref{eq:system1D}. However, the variance $\lambda^2$ plays the role of a bifurcation parameter since for $\lambda^2 \to 0^+$ we recover the Bogdanov-Takens-type equilibria and thus $\lambda^2$ changes the stability of the equilibrium point. In view of this consideration we wish to avoid $\lambda^2$ going to zero and, furthermore, we can apply a control on it in order to guarantee that the unstable equilibrium $\tilde{F}_1^-$ is always negative and thus not admissible.  More precisely, the standard deviation should satisfy
$$
\lambda > \frac{y^2\sqrt{2}}{G}.
$$
Then, the solutions of~\eqref{eq:systemNoise1D} are given by
$$
m(t) = \frac{y}{G}, \quad  E(t) = \frac{\left(\tanh(\sqrt{2\lambda^2} G C + \sqrt{2\lambda^2} G t \right) \sqrt{2\lambda^2} G + 2 y^2}{2 G^2}
$$
and
\begin{align*}
m(t) &= \frac{\pm e^{2 \sqrt{2\lambda^2} G t} C_1 y \mp C_2 y + \sqrt{\sqrt{2\lambda^2} C_1 e^{3\sqrt{2\lambda^2} G t} - C_2 \sqrt{2\lambda^2} e^{2\sqrt{2\lambda^2} G t}}}{(C_1 e^{2\sqrt{2\lambda^2} G t} - C_2) G}, \\
E(t) &= \frac{m(t)^3 G^2 - m(t)^2 G y - \frac{\mathrm{d}}{\mathrm{d}t} m(t)}{m(t) G^2 - G y}
\end{align*}
with $C, C_1, C_2 \in \R$ constants uniquely prescribed by the initial conditions. We observe that, in the large time behavior, $m \to \frac{y}{G}$ unconditionally, as in the case of~\eqref{eq:system1D}. Instead, the large time behavior of $E$ is changed and shifted by a quantity depending on $\lambda^2$ which avoids the possibility of having a decay in the energy which drives $m$ away from the expected equilibrium value. In Figure~\ref{fig:phaseplaneNoise} we show the nullclines and the complete vector field for the case $(y,G)=(2,1)$ (left panel) and the behavior around the stable equilibrium point $\tilde{F}_1^+ = (2,8)$ (right panel).

\begin{remark}
	Lemma~\ref{lem:schillingsstuart} shows that the collapse of the ensembles towards their mean slows down linearly as the number of the ensemble increases. The kinetic equation~\eqref{eq:kineticFP} holds in the limit of a large ensemble size and the energy $\vec{E}$ gives information on the concentration of the distribution $f$ of the control $\vec{u}$ around its mean. The previous analysis shows also that, in fact, the result of Lemma~\ref{lem:schillingsstuart} holds at the kinetic level since $\vec{E}$ does not decay to zero as $t\to\infty$.
\end{remark}

\section{Numerical simulation results} \label{sec:numericalTest}

The simulations are performed by using a standard Monte Carlo approach~\cite{Caflisch1998} to solve the kinetic equation~\eqref{eq:kineticFP}. More precisely, we use a simple modification of the mean-field interaction algorithm given in~\cite{AlbiPareschi2013} which is a direct simulation Monte Carlo method based on the mean-field microscopic dynamics described by~\eqref{eq:microInteraction} giving rise to the corresponding kinetic equation~\eqref{eq:kineticFP}. For further details on the method we refer to ~\cite{BabovskyNeunzert1986,FornasierEtAl2011,Lemou1998,MouhotPareschi2006,PareschiRusso1999,PareschiToscaniBOOK}. 

The algorithmic details are as follows. In each example we consider a sampling of $J$ controls $\{\vec{u}^j\}_{j=1}^J$ from the prior or initial distribution $f_0(\vec{u})$. Then, each sample is updated according to the mean-field microscopic rule~\eqref{eq:microInteraction} by selecting $M \leq J$ interacting particles uniformly distributed without repetition. The parameter $\epsilon$ in~\eqref{eq:microInteraction} is closely related with the concept of a time step and it is taken such that  stability of the discrete method is guaranteed.~\cite{AlbiPareschi2013}. In particular,  for the kinetic model~\eqref{eq:kineticFromEnKF} we require that
\begin{equation} \label{eq:stabilitycondition}
	\epsilon \leq \frac{1}{\max_{i}\left(|(\Re(\mu_i)|\right)}
\end{equation}
where the $\mu_i$'s are the eigenvalues of $\vecsym{\C}(t) G^T \vecsym{\Gamma} G$, cf. Remark~\ref{rem:probabilisticMicro}. As we observe that $\vecsym{\C}(t) G^T \vecsym{\Gamma} G$ is characterized by large spectral radius at initial time that reduces over time, we chose an adaptive computation of $\epsilon$ by recomputing it at each iteration.

As already pointed out in Section~\ref{sec:stabilization}, the microscopic interactions~\eqref{eq:microInteraction} are closely related to a time discretization of the gradient descent equation~\eqref{eq:gradientEq}. However, a deterministic numerical method for~\eqref{eq:gradientEq} requires $O(J^2)$ operations due to the direct evaluation of the sum for $J$ ensembles. The numerical discretization of the kinetic equation by means of a Monte Carlo approach allows to compute the microscopic dynamics with a cost directly proportional to the number $J$ of ensembles.

Information on the simulation results is presented in  the following norms:
\begin{equation} \label{eq:normErr}
\begin{aligned}
	v &= \frac{1}{J} \sum_{j=1}^J \| \vec{v}^j \|_2^2, \quad r = \frac{1}{J} \sum_{j=1}^J \| \vec{r}^j \|_2^2\\
	V &= \frac{1}{J} \sum_{j=1}^J | \vec{V}_{jj} |^2, \quad R = \frac{1}{J} \sum_{j=1}^J | \vec{R}_{jj} |^2
\end{aligned}
\end{equation}
which are computed at each iteration and where
\begin{equation} \label{eq:singleErr}
\begin{aligned}
	\vec{v}^j &= \vec{u}^j - \overline{\vec{u}}, \quad \vec{r}^j = \vec{u}^j - \vec{u}^\dagger \\
	\vec{V}_{ij} &= \left\langle G \vec{v}^i, G \vec{v}^j \right\rangle_{\vecsym{\Gamma}^{-1}}, \quad \vec{R}_{ij} = \left\langle G \vec{r}^i, G \vec{r}^j \right\rangle_{\vecsym{\Gamma}^{-1}}.
\end{aligned}	
\end{equation}
The quantity $\vec{v}^j$ measures the deviation of the $j$-th sample from the mean $\overline{\vec{u}}$ of the approximated distribution by the samples and $\vec{r}^j$ measures the deviation of the $j$-th sample from the truth solution $\vec{u}^\dagger$. The quantities $\vec{V}$ and $\vec{R}$ give information on the deviation of $\vec{v}^j$ and $\vec{r}^j$ under application of the model $G$.

Another additional important quantities is given by the misfit which allows to measure the quality of the solution at each iteration. The misfit for the $j$-th sample is defined as
\begin{equation} \label{eq:singleMisfit}
	\vecsym{\vartheta}^j = G \vec{r}^j - \vecsym{\eta}.
\end{equation}
By using~\eqref{eq:singleMisfit} we finally look at
\begin{equation} \label{eq:misfit}
	\vartheta = \frac{1}{J} \sum_{j=1}^J \| \vecsym{\vartheta}^j \|_{\vecsym{\Gamma^{-1}}}^2.
\end{equation}
Driving this quantity to zero leads to over-fitting of the solution. For this reason, usually it is suitable introducing a stopping criterion which avoids this effect. In the following we will consider the discrepancy principle which check and stop the simulation when the condition $\vartheta \leq \| \vecsym{\eta} \|_2^2$ is satisfied.

The algorithm employed in the experiments is summarized by the steps described in Algorithm~\ref{alg:meanfield}.

\begin{algorithm}[t!]
	\begin{algorithmic}[1]
		\STATE Given $J$ samples $\vec{u}^{j,0}$, with $j=1,\dots,J$ computed from the initial distribution $f_0(\vec{u})$ and $M\leq J$;
		\STATE set $n=0$, $t^0=0$ and a final time $T_\text{fin}$;
		\WHILE{$n\leq n_{tot}$}
		\STATE compute the misfit $\vartheta$ as in~\eqref{eq:misfit};
		\IF{$\vartheta \leq \| \vecsym{\eta} \|_2^2$}
		\BREAK
		\ELSE
		\STATE compute $\epsilon = \frac{1}{\max_{i}\left(|(\Re(\mu_i)|\right)}$ $\mu_i$'s are the eigenvalues of $\vecsym{\C}(t) G^T \vecsym{\Gamma} G$;
		\IF{$t^n+\epsilon>T_\text{fin}$}
		\STATE set $\epsilon = T_\text{fin}-t^n$;
		\ENDIF
		\IF{$t^{n}\geq T_\text{fin}$}
		\BREAK
		\ENDIF
		\FOR{$j=1$ \TO $J$}
		\STATE sample $M$ data $j_1,\dots,j_M$ uniformly without repetition among all data;
		\STATE compute
		$$
		\vec{m}_M^n = \frac{1}{M} \sum_{k=1}^M \vec{u}^{j_k,n}, \quad \vec{E}_M^n = \frac{1}{M} \sum_{k=1}^M \vec{u}^{j_k,n} \otimes \vec{u}^{j_k,n};
		$$
		\STATE sample $\vecsym{\xi}$ from a zero mean distribution $\theta(\vecsym{\xi})$ having given covariance matrix $\vecsym{\Lambda}$;
		\STATE compute the data change
		$$
		\vec{u}^{j,n+1} = \vec{u}^{j,n} - \epsilon (\vec{E}_M^n - \vec{m}_M^n \otimes \vec{m}_M^n) \nabla_\vec{u} \Phi(\vec{u}^{j,n},\vec{y}) + \sqrt{\epsilon} \, \vecsym{\xi};
		$$
		\ENDFOR
		\ENDIF
		\STATE set $n=n+1$ and $t^{n+1}=t^n+\epsilon$.
		\ENDWHILE
	\end{algorithmic}
	\caption{Mean-Field Ensemble Kalman Filter Algorithm.}
	\label{alg:meanfield}
\end{algorithm}

\subsection{Linear elliptic problem} \label{sec:elliptic}

A test proposed e.g. in~\cite{iglesiaslawstuart2013,schillingsstuart2017}, is the ill-posed inverse problem of finding the force function of an elliptic equation in one spatial dimension assuming that  noisy observation of the solution to the problem are available. This  problem is widely used since is explicitly solvable due to the linearity of the model.

The problem is prescribed by the following one dimensional elliptic equation
$$
	-\frac{\mathrm{d}^2}{\mathrm{d}x^2} p(x) + p(x) = u(x), \quad x\in[0,\pi]
$$
endowed with boundary conditions $p(0) = p(\pi) = 0$. The linear model is thus defined as
$$
	A = \left( -\frac{\mathrm{d}^2}{\mathrm{d}x^2} + 1 \right)^{-1}
$$
which can be discretized, for instance, by a finite difference method or by the explicit solution 
$$
	p(x) = A \, u(x) = \exp(x) \left( C_1 - \frac12 \int_0^x \exp(y) u(y) \mathrm{d}y \right) + \exp(-x) \left( C_2 + \frac12 \int_0^x \exp(-y) u(y) \mathrm{d}y \right)
$$
where the constants $C_1$ and $C_2$ can be uniquely determined by the boundary conditions. We assign a continuous  control $u(x)$ and then  introduce a uniform mesh consisting of $d=K=2^8$ equidistant points in the interval $[0,\pi]$. Let $\vec{u}^\dagger\in\R^d$ be the vector of the evaluations of the control function $u(x)$ on the mesh. We simulate noisy observations $\vec{y}\in\R^K$ as
$$
	\vec{y} = \vec{p} + \vecsym{\eta} = G \vec{u}^\dagger + \vecsym{\eta},
$$
where $G$ is the finite difference discretization of the continuous operator $A$. For simplicity we assume that $\vecsym{\eta}$ is a Gaussian white noise, more precisely $\vecsym{\eta}\sim \mathcal{N}(0,\gamma^2 \vec{I})$ with $\gamma \in \R^+$ and $\vec{I} \in \R^{d\times d}$ is the identity matrix. We are interested in recovering the control $\vec{u}^\dagger \in \R^d$ from the noisy observations $\vec{y}\in\R^K$ only.

The initial ensemble of particles is sampled by an initial distribution $f_0(\vec{u}) = \mathcal{N}(0,\vec{C}_0)$. The choice of $f_0(\vec{u})$ is related to the choice of the prior distribution in Bayesian problems. In this case $f_0(\vec{u})$ represents a Brownian bridge as in~\cite{schillingsstuart2017}


\paragraph{Test case 1.}
\begin{figure}[t!]
	\centering
	\includegraphics[width=0.49\textwidth]{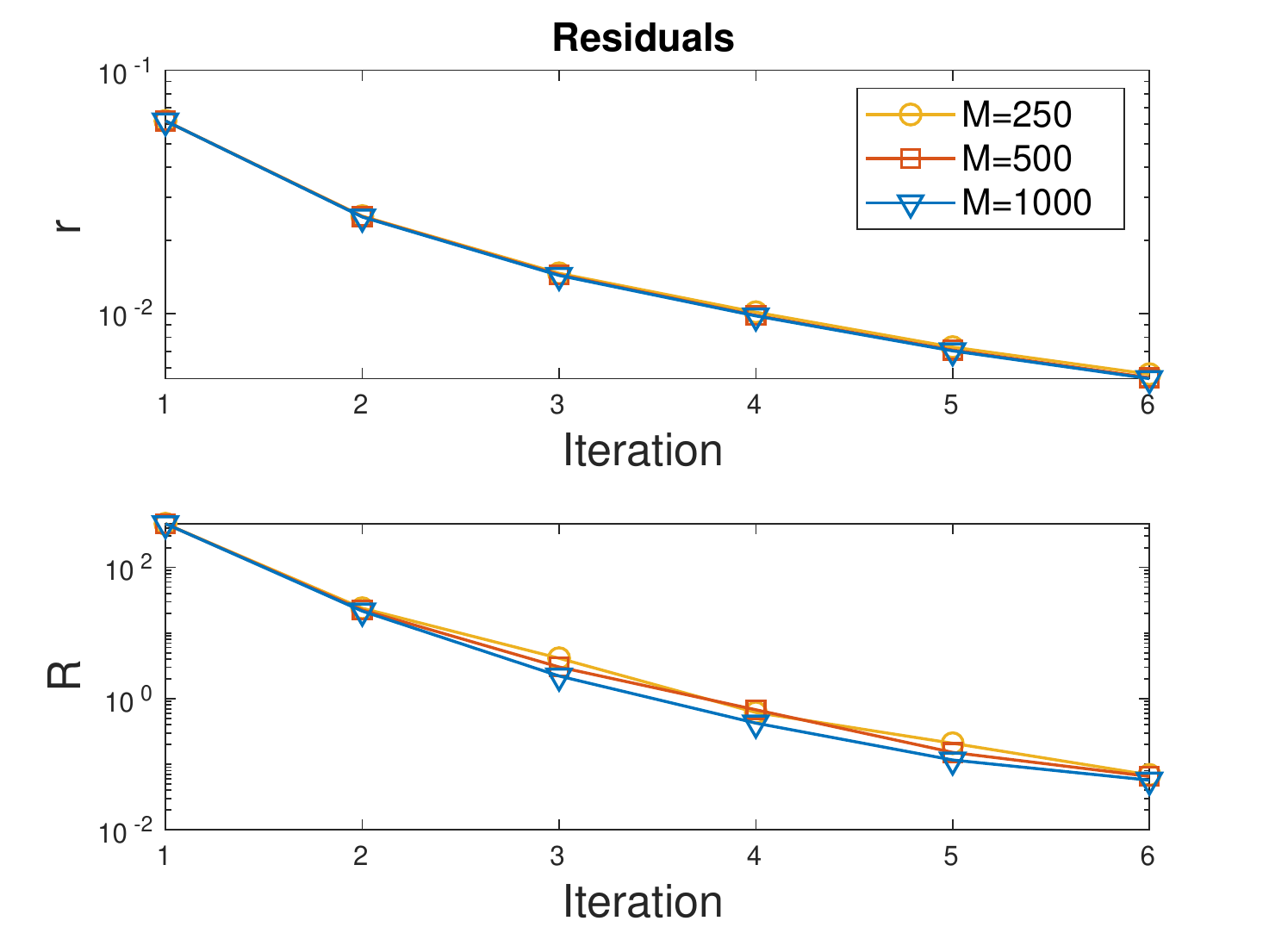}
	\includegraphics[width=0.49\textwidth]{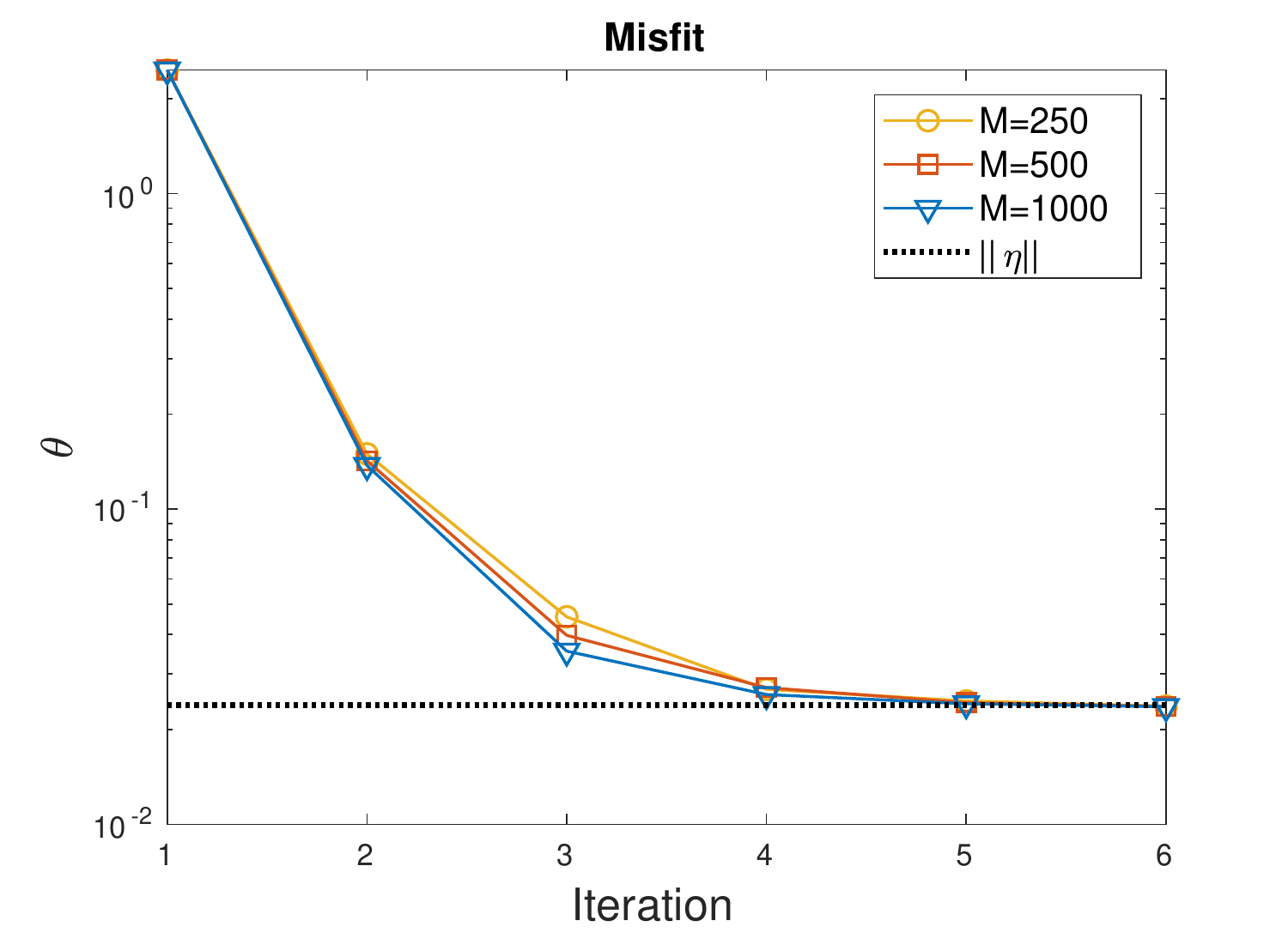}
	\\
	\includegraphics[width=0.32\textwidth]{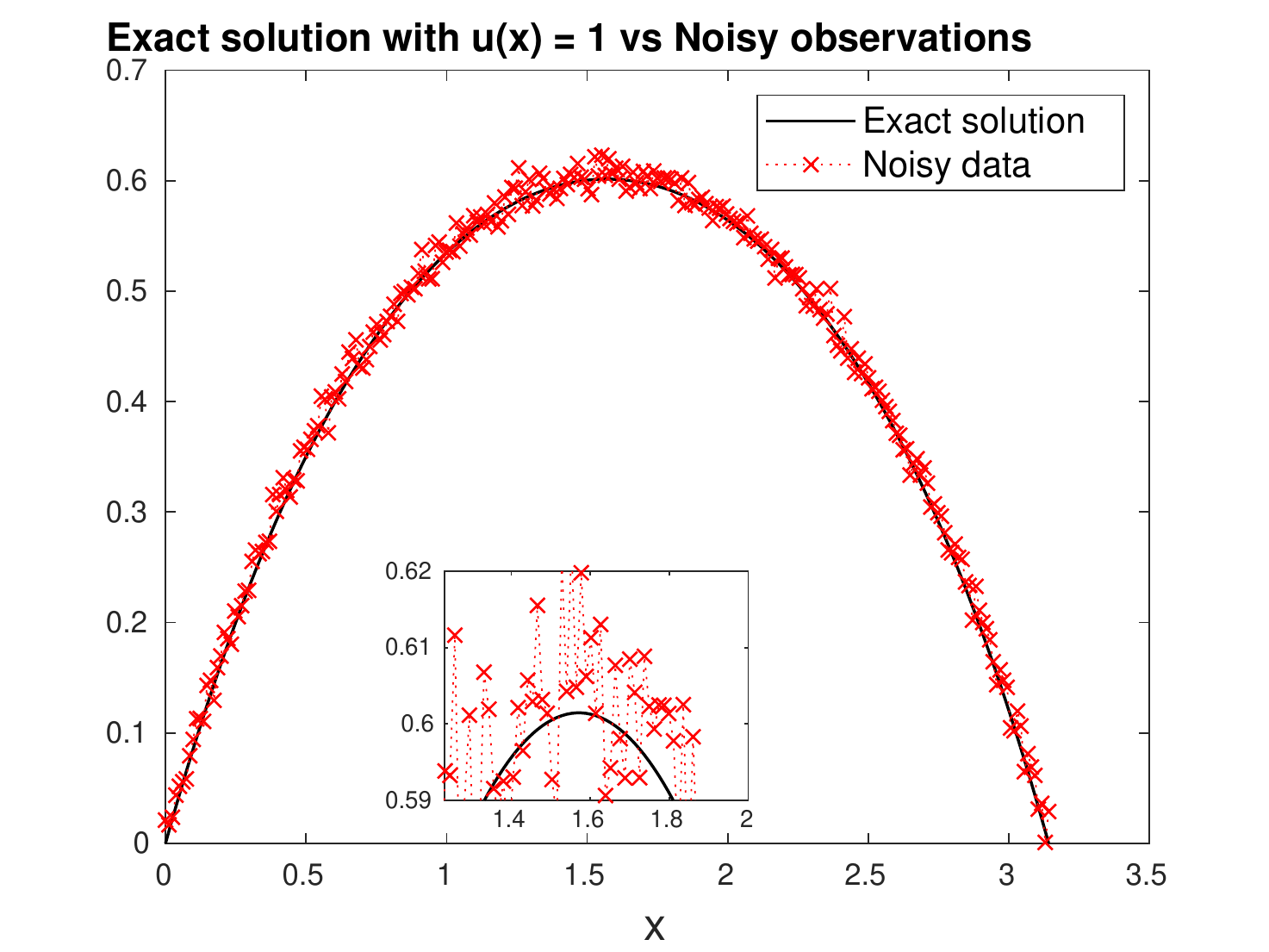}
	\includegraphics[width=0.32\textwidth]{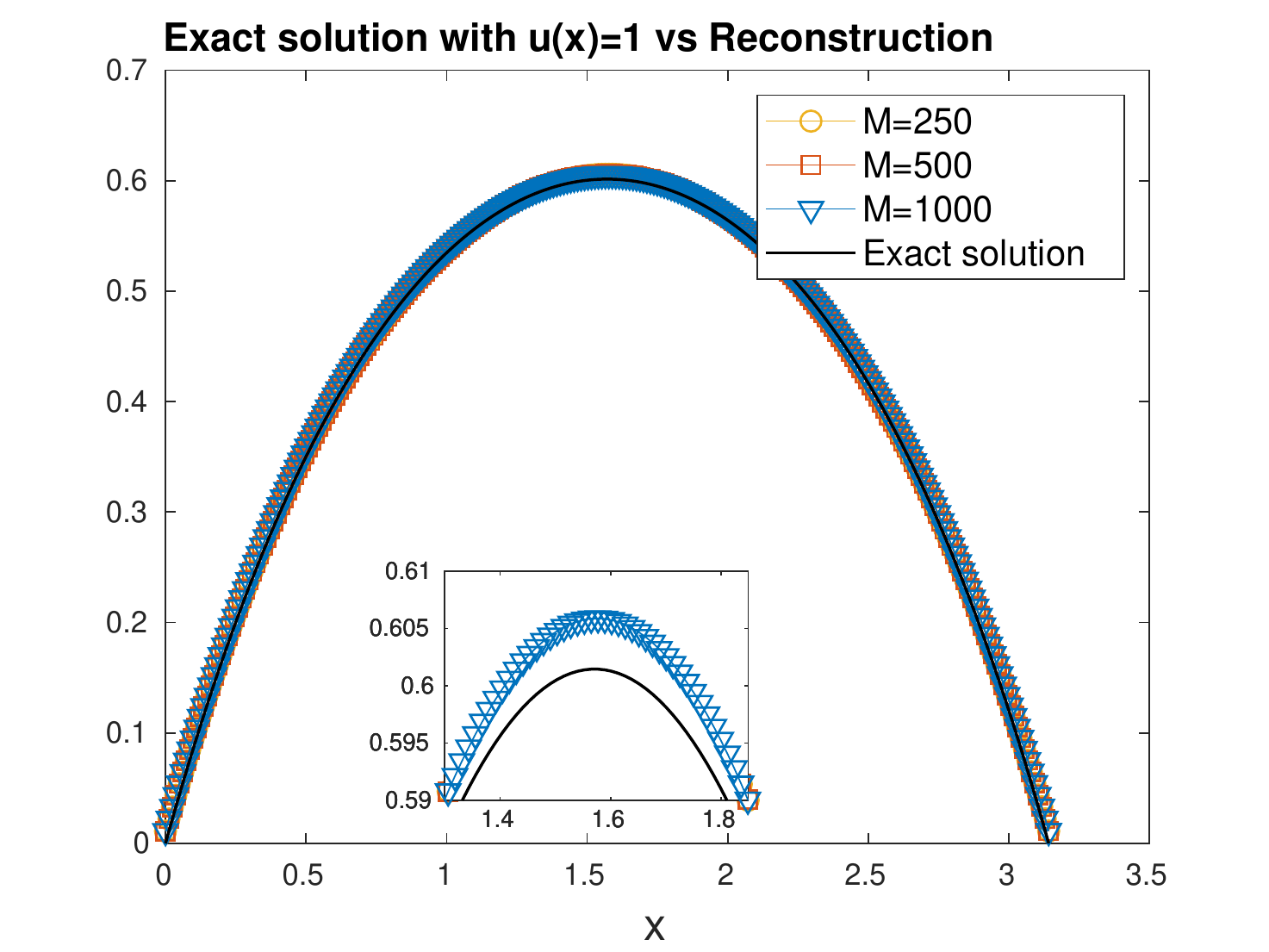}
	\includegraphics[width=0.32\textwidth]{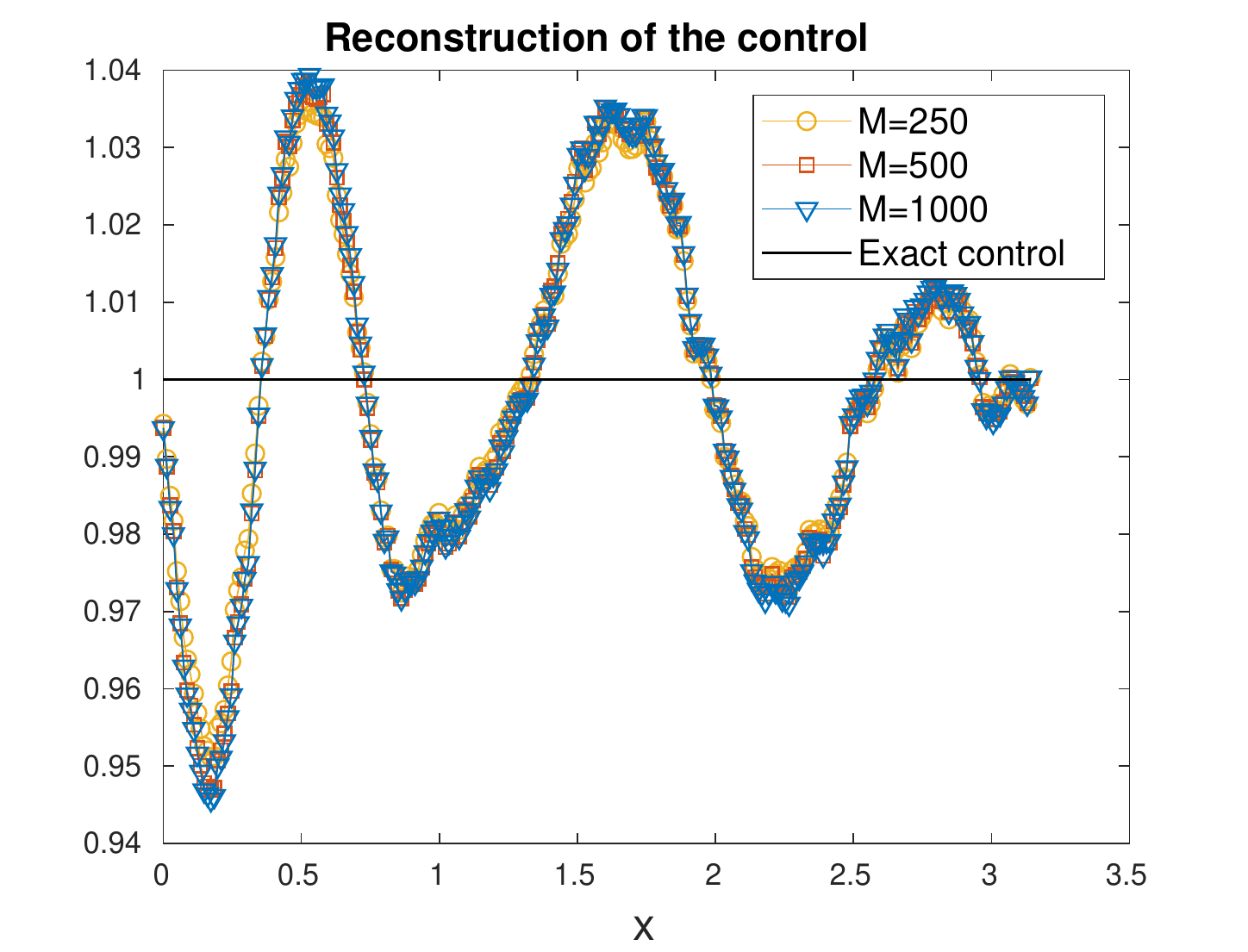}
	\caption{Elliptic problem - Test case 1 with $\gamma=0.01$. Top row: plots of the residual $r$, the projected residual $R$ and the misfit $\vartheta$ for $M=250,500,1000$. Bottom row: plots of the noisy data, the reconstruction of $p(x)$ and the reconstruction of the control $u(x)$ at final iteration for $M=250,500,1000$.\label{fig:test1gamma001}}
\end{figure}

\begin{figure}[t!]
	\centering
	\includegraphics[width=0.49\textwidth]{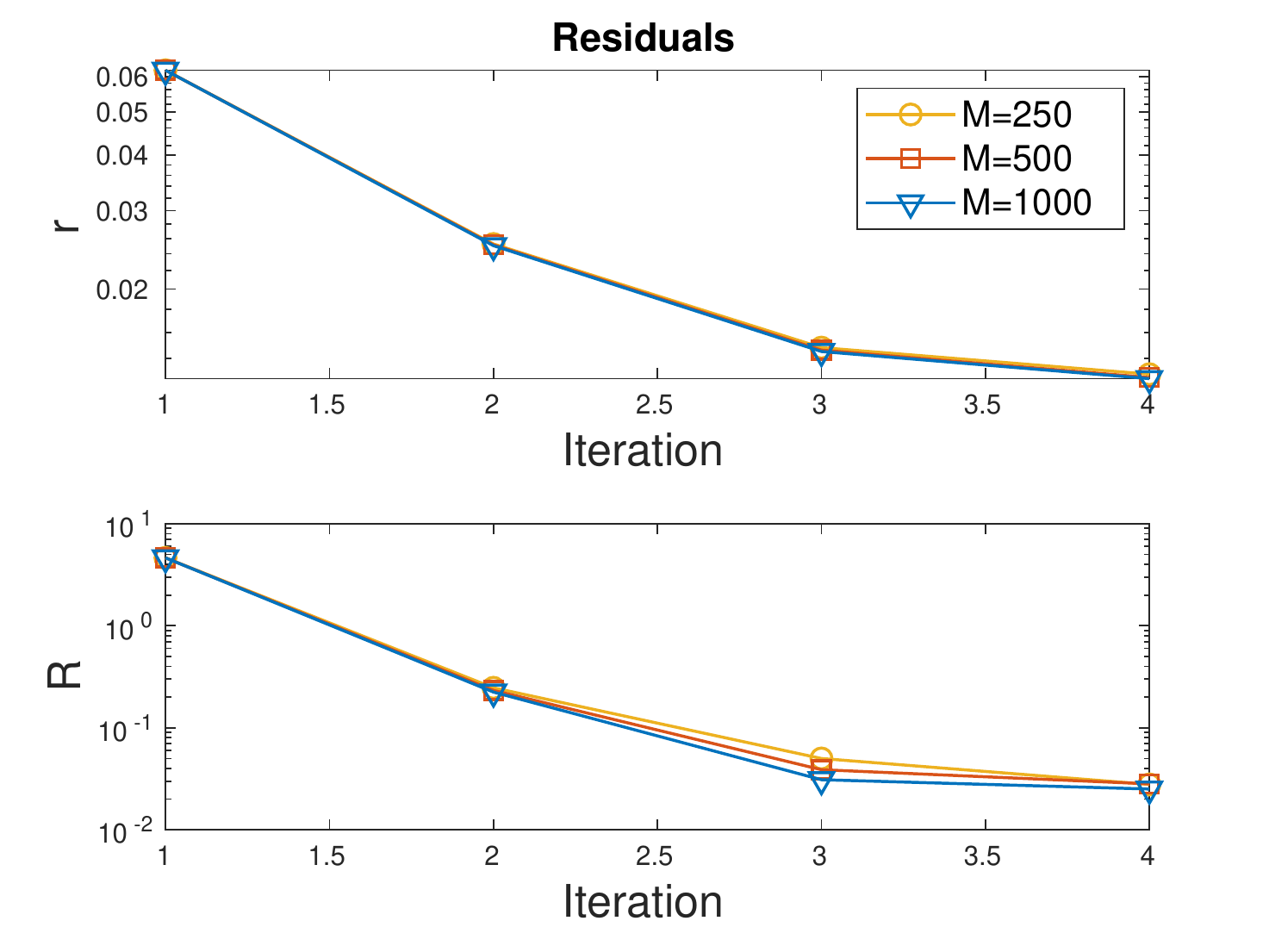}
	\includegraphics[width=0.49\textwidth]{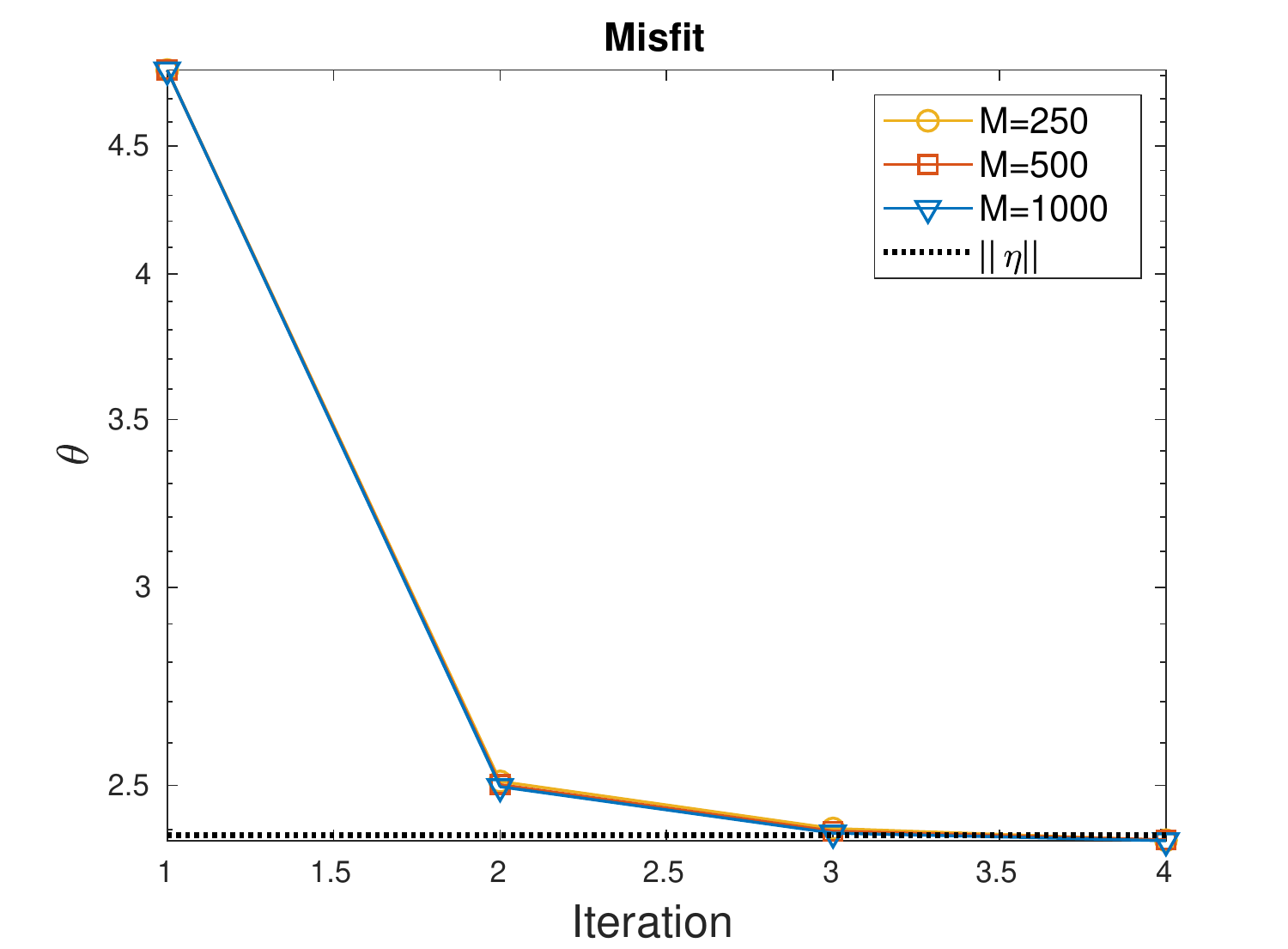}
	\\
	\includegraphics[width=0.32\textwidth]{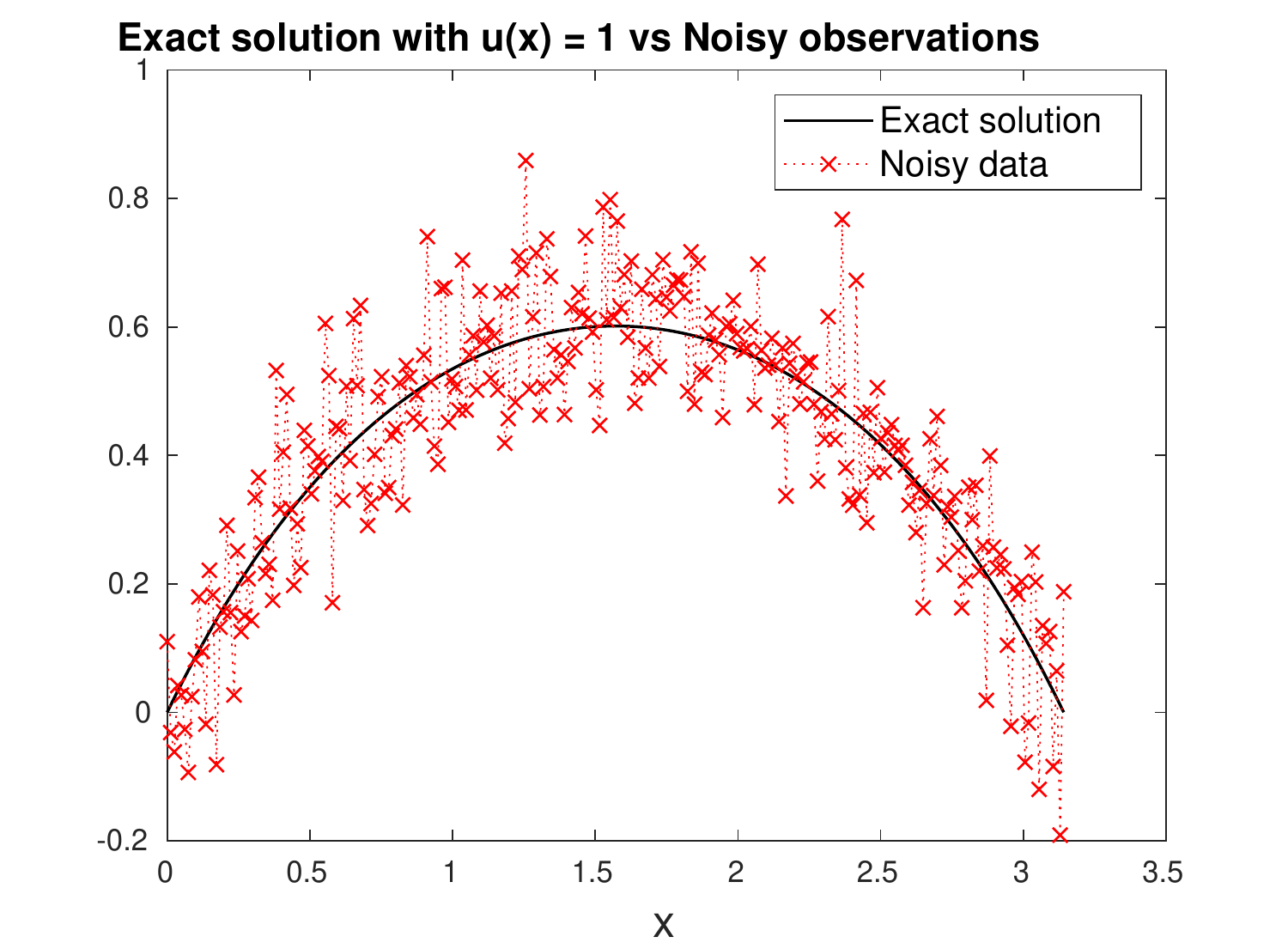}
	\includegraphics[width=0.32\textwidth]{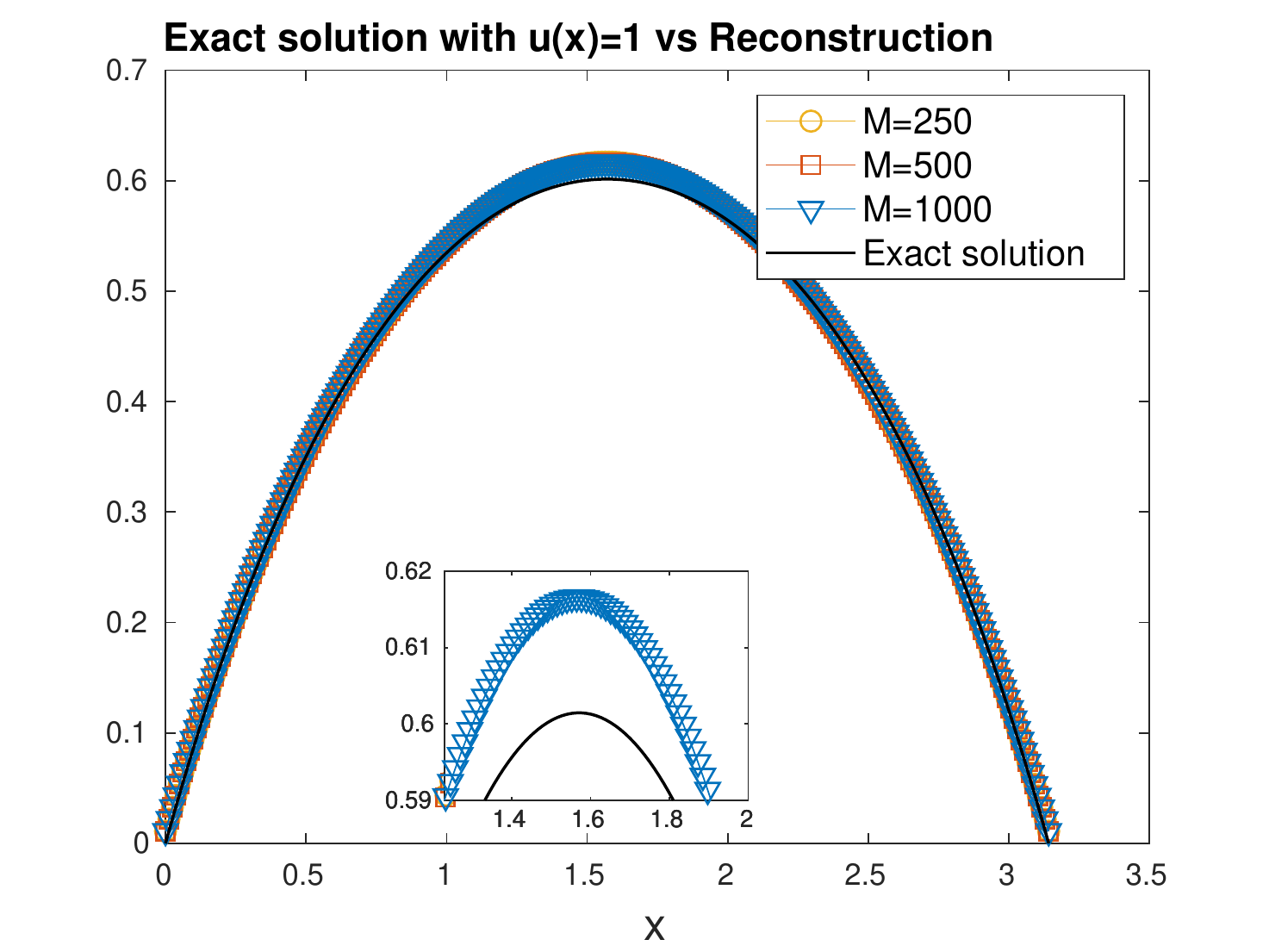}
	\includegraphics[width=0.32\textwidth]{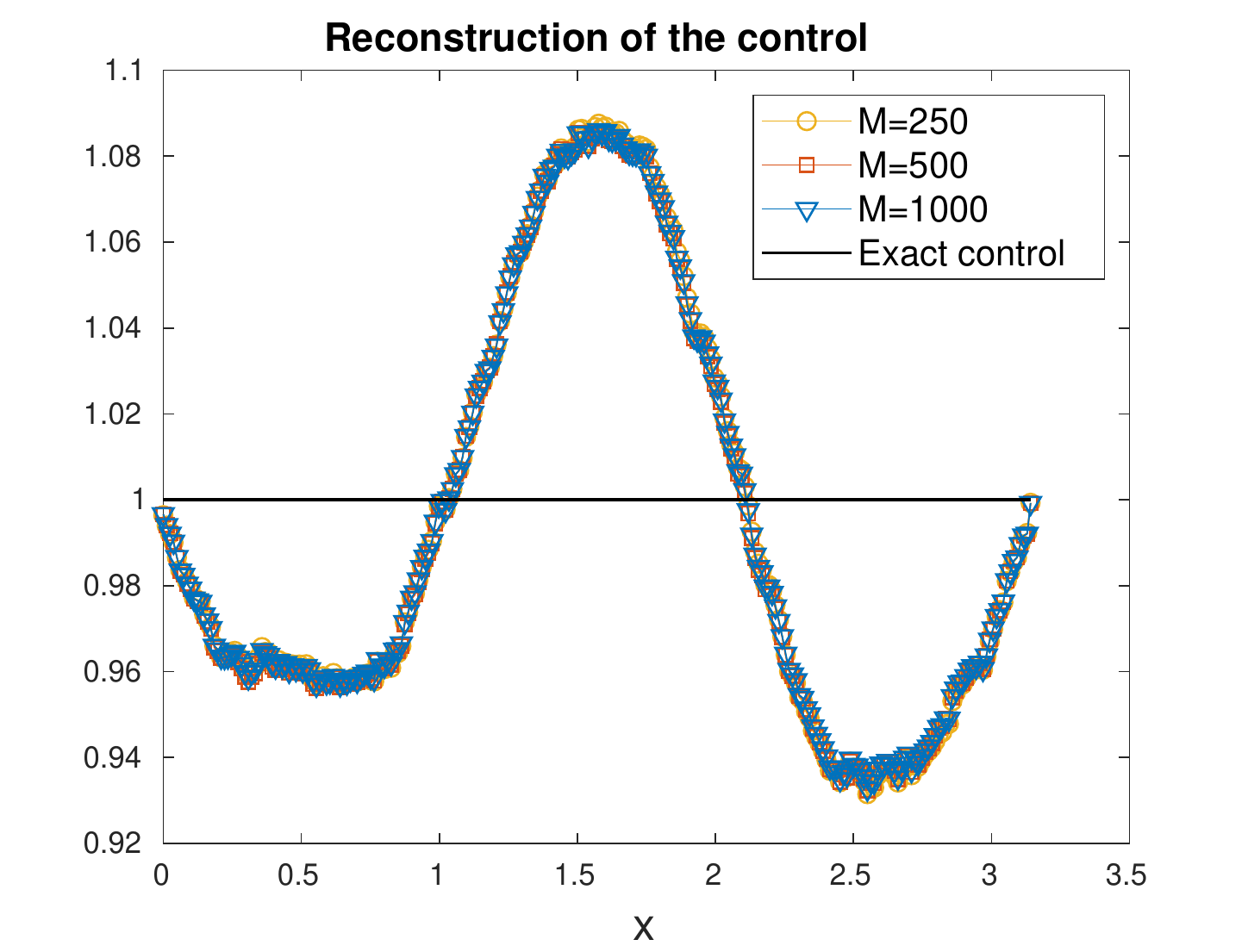}
	\caption{Elliptic problem - Test case 1 with $\gamma=0.1$. Top row: plots of the residual $r$, the projected residual $R$ and the misfit $\vartheta$ for $M=250,500,1000$. Bottom row: plots of the noisy data, the reconstruction of $p(x)$ and the reconstruction of the control $u(x)$ at final iteration for $M=250,500,1000$.\label{fig:test1gamma01}}
\end{figure}

\begin{figure}[t!]
	\centering
	\includegraphics[width=0.49\textwidth]{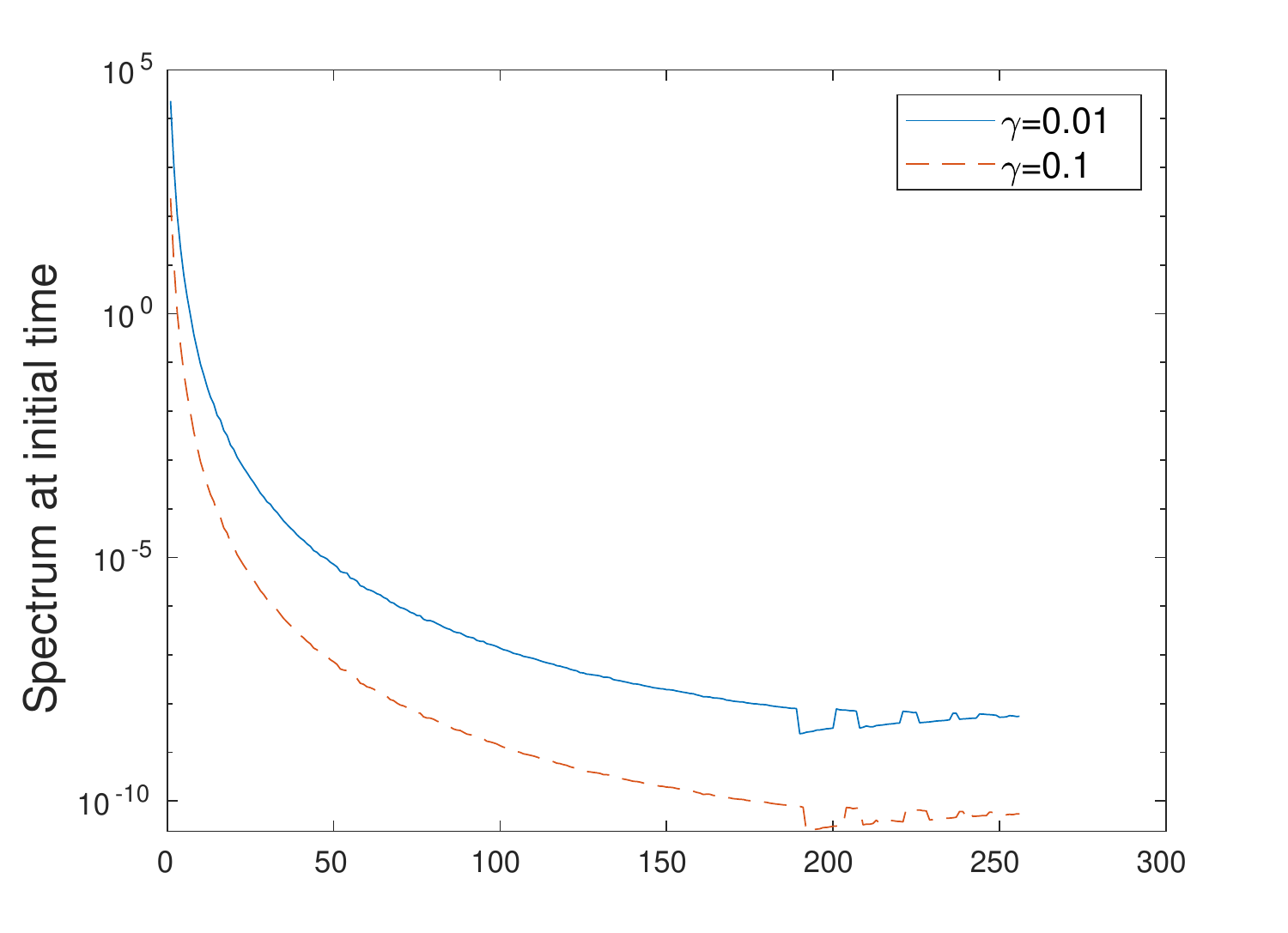}
	\includegraphics[width=0.49\textwidth]{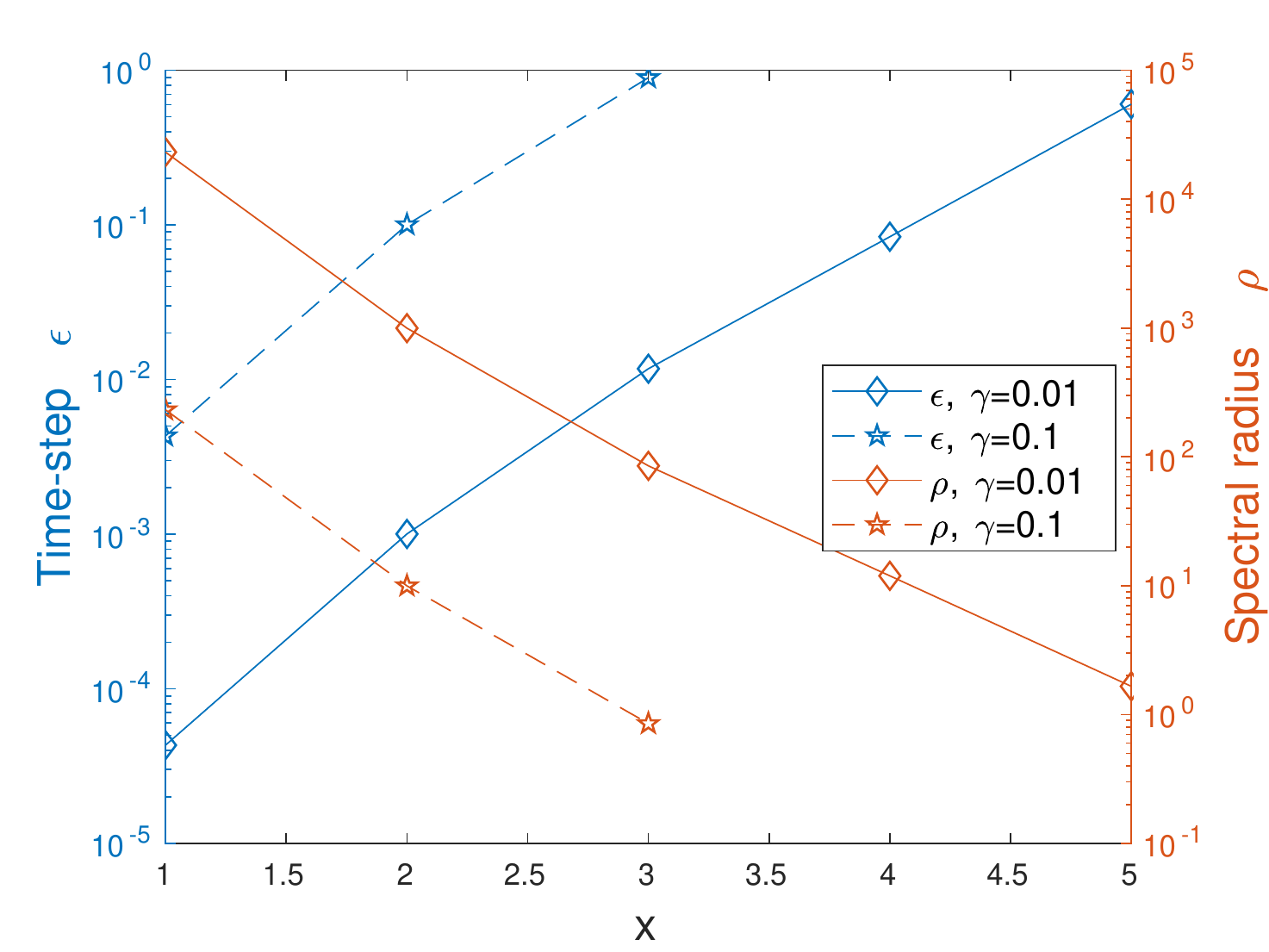}
	\caption{Elliptic problem - Test case 1. Left: spectrum of $\vecsym{\C}(t) G^T \vecsym{\Gamma} G$ for the initial data with $\gamma=0.01$ and $\gamma=0.1$. Right: adaptive $\epsilon$ and spectral radius of $\vecsym{\C}(t) G^T \vecsym{\Gamma} G$ over iterations with $\gamma=0.01$ and $\gamma=0.1$.\label{fig:test1spectral}}
\end{figure}

Let us consider $u(x) = 1$, $\forall\,x\in[0,\pi]$. We solve the inverse problem by the proposed method for  different values $M$ of the interacting samples. We observe that taking $M < J$ does not strongly influence the results of the simulation. But, $M < J$ allows to have a computational gain.

We allow for two values of the noise level $\gamma=0.01$, see Figure~\ref{fig:test1gamma001}, and $\gamma=0.1$, see Figure~\ref{fig:test1gamma01}. In both figures, the top panels show the residual (left) and misfit (right) decrease over the number of iterations. Due to the discrepancy principle, the simulation is automatically stopped when the misfit reaches $\|\vecsym{\eta}\|$. The final residual values are obviously larger in the case of $\gamma=0.1$ due to the larger noise level present in the initial observations. The bottom panels show, form left to right, the initial noisy data which are spread around the exact solution $p(x)$ of the problem, the reconstruction of $p(x)$ and the reconstruction of the control $u(x)$ by using the mean of the samples as estimator of the solution. We observe that the different values of $M$ does not give significantly different results.

In the left panel of Figure~\ref{fig:test1spectral} we show the spectrum of $\vecsym{\C}(t) G^T \vecsym{\Gamma} G$ at initial time for $\gamma=0.01$ and $\gamma=0.1$. We observe that the ratio between the largest and smaller eigenvalues is very large, reflecting the ill-posedness of the problem and the need of using a small $\epsilon$ in~\eqref{eq:probabInteraction} in order to guarantee stability. However, we consider an adaptive $\epsilon$ since the spectral radius is observed to decrease quickly over iterations. See the red lines in the right panel of Figure~\ref{fig:test1spectral}, where, instead, the blue lines show the corresponding values of $\epsilon$ which avoid the lack of stability.

\paragraph{Test case 2.}

\begin{figure}[t!]
	\centering
	\includegraphics[width=0.49\textwidth]{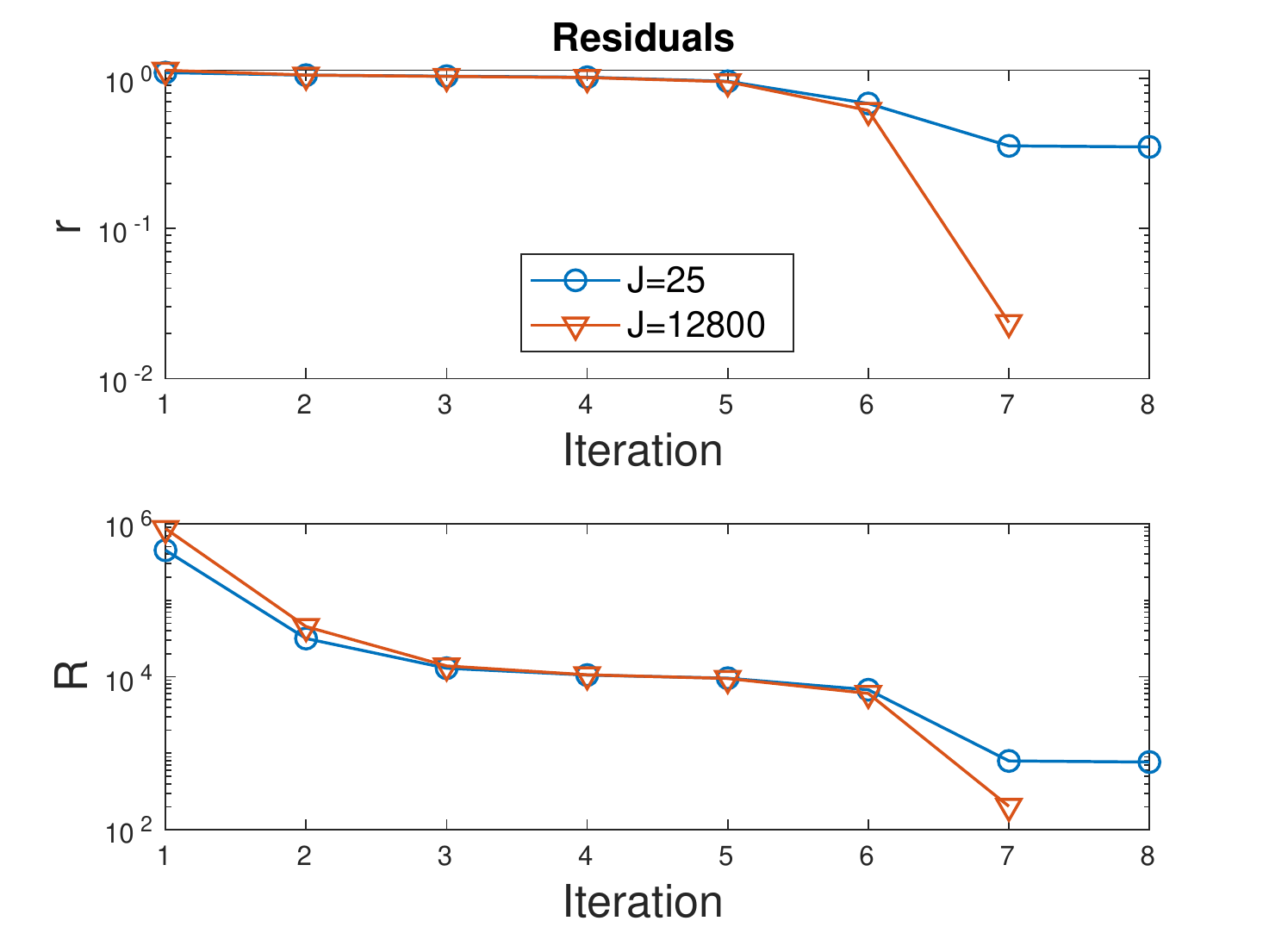}
	\includegraphics[width=0.49\textwidth]{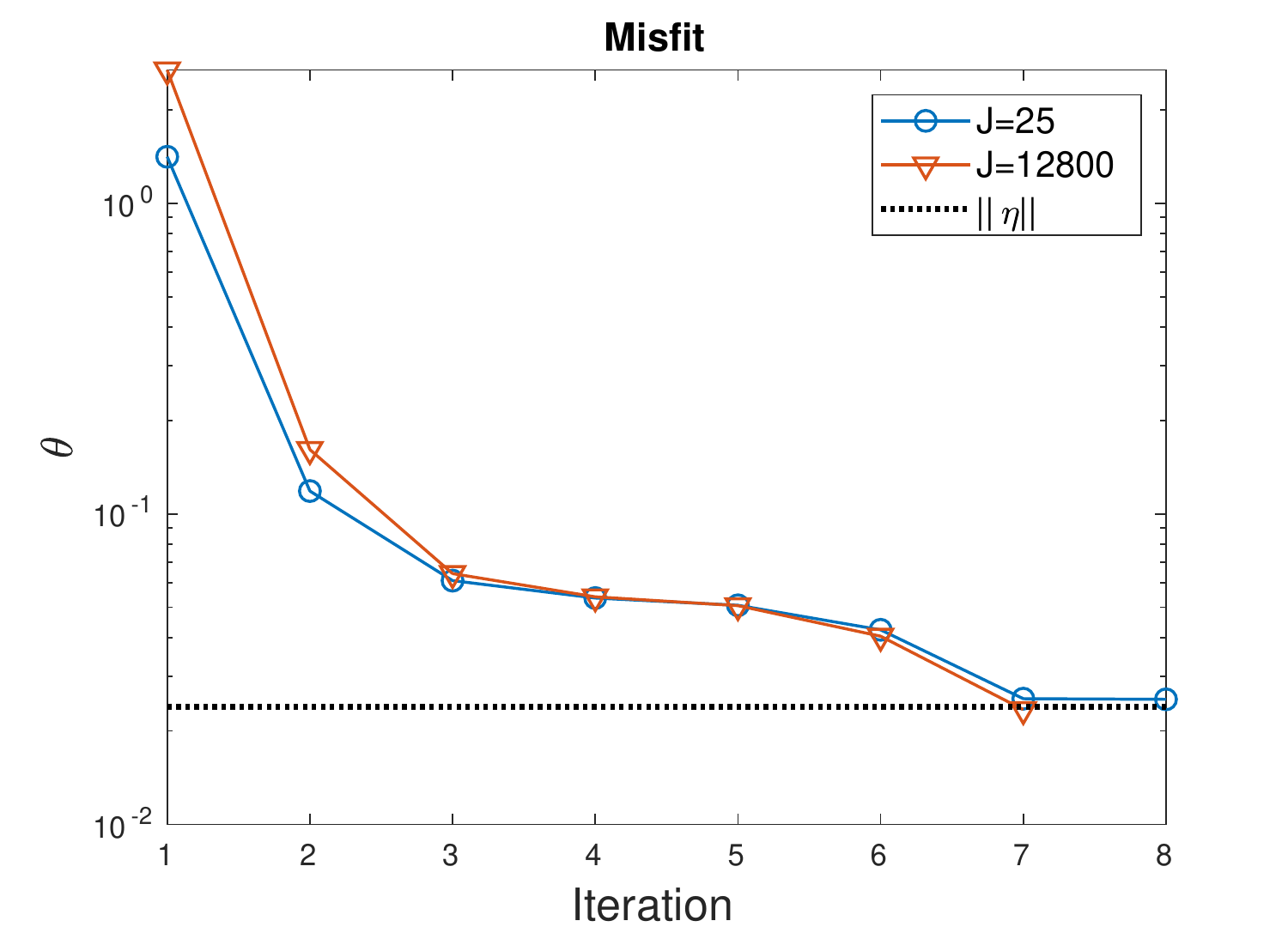}
	\\
	\includegraphics[width=0.49\textwidth]{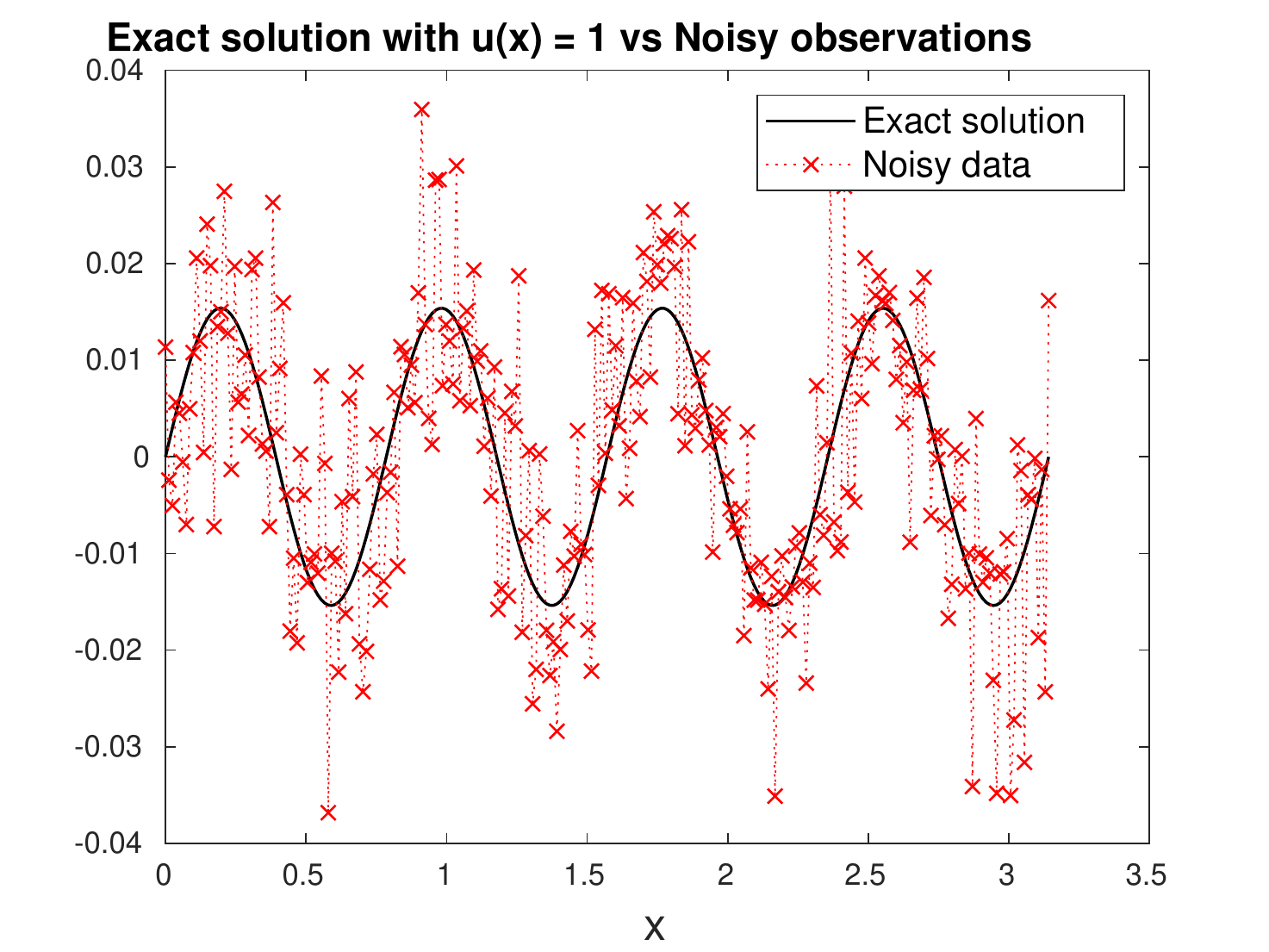}
	\includegraphics[width=0.49\textwidth]{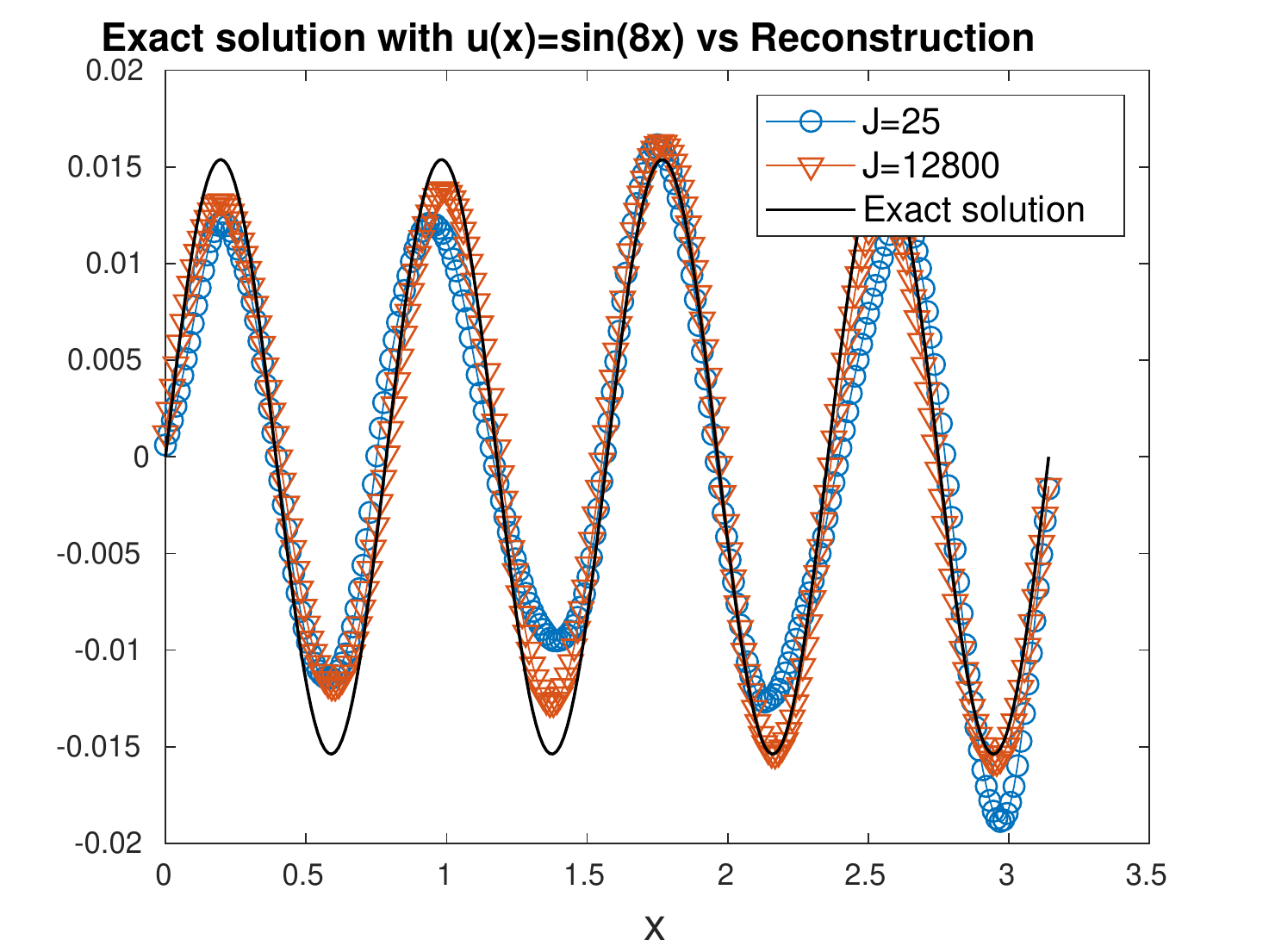}
	\\
	\includegraphics[width=0.49\textwidth]{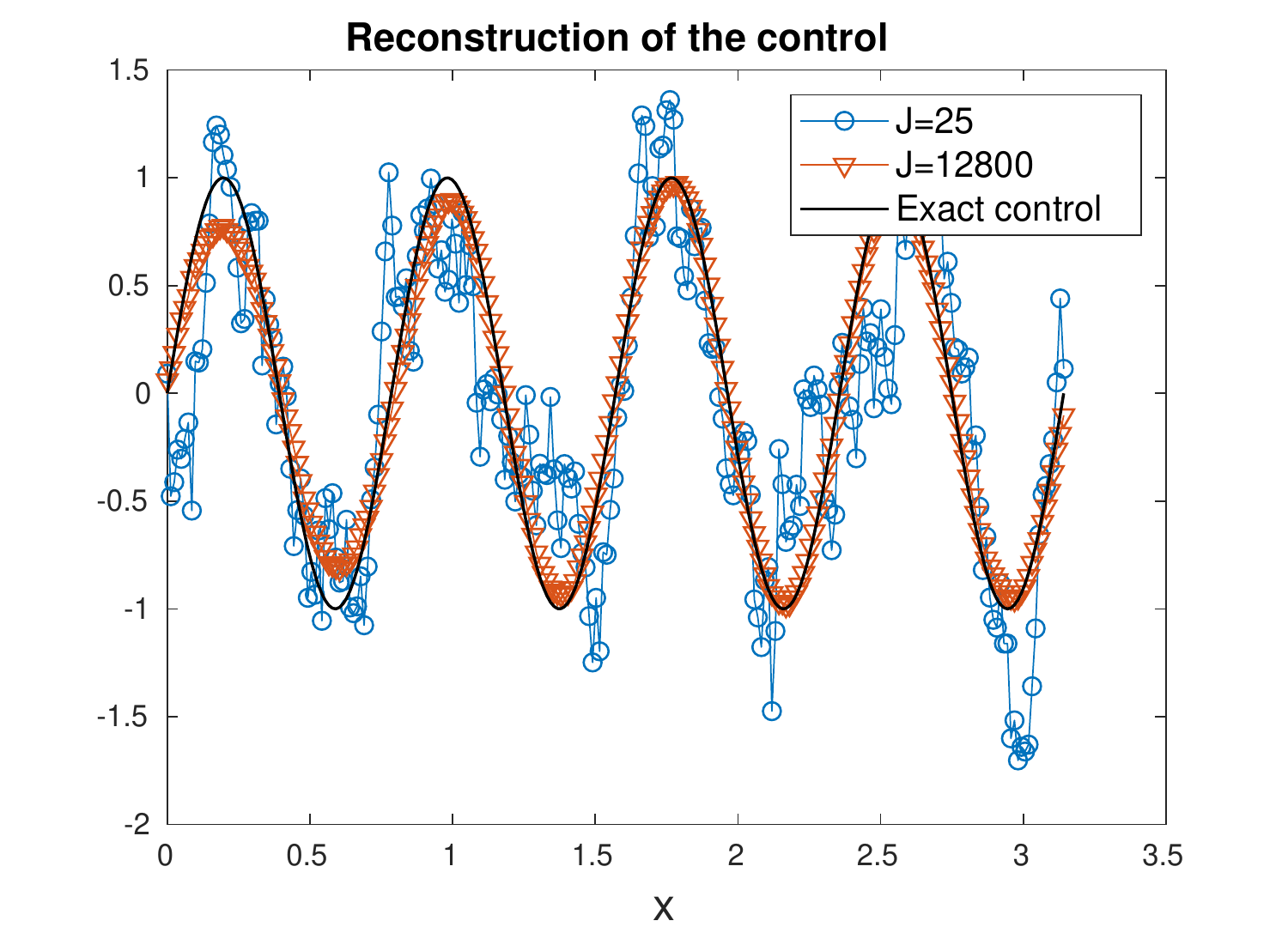}
	\includegraphics[width=0.49\textwidth]{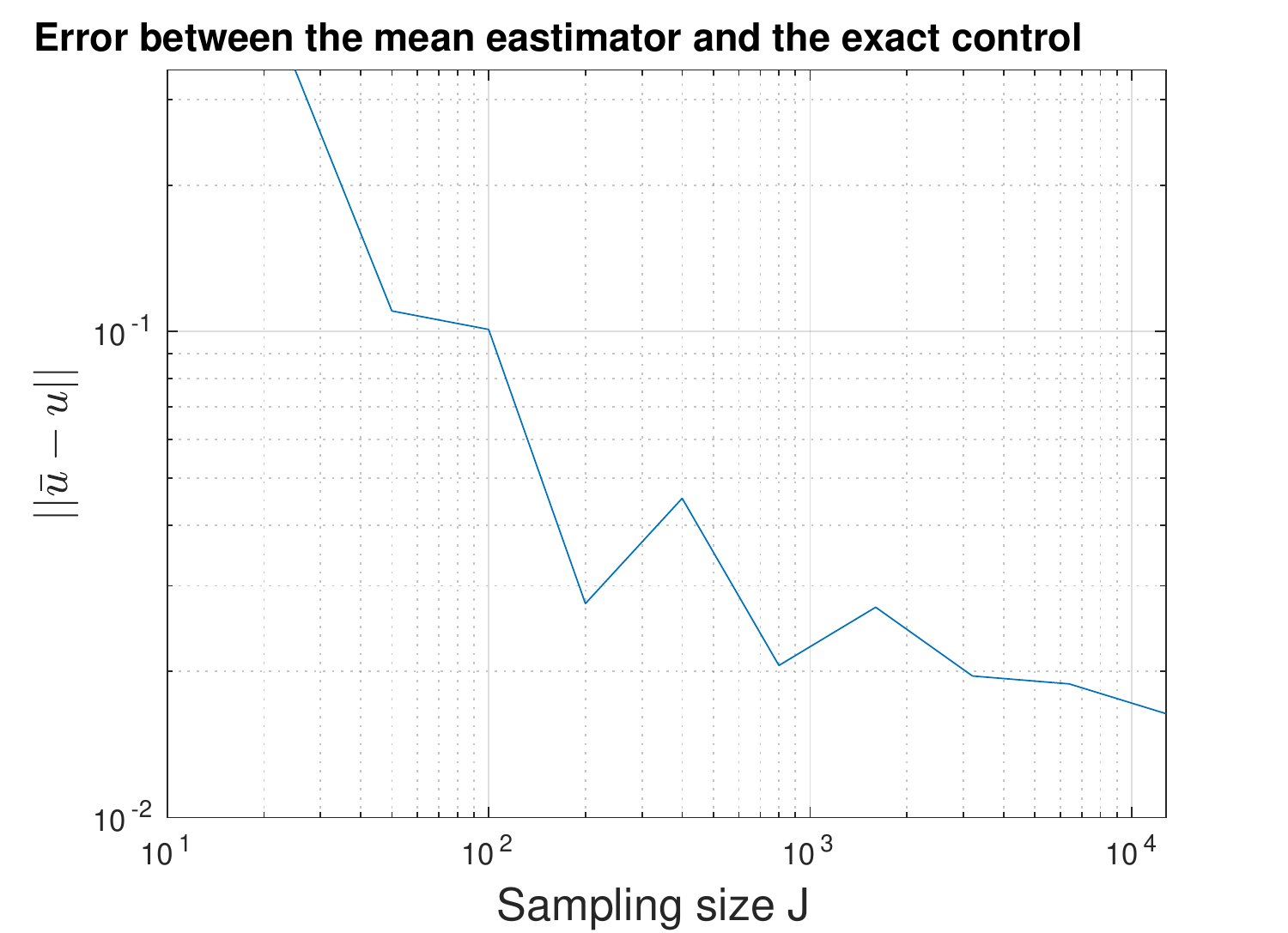}
	\caption{Elliptic problem - Test case 2 with $\gamma=0.01$. Top row: plots of the residual $r$, the projected residual $R$ and the misfit $\vartheta$ for $J=25,25\cdot2^9$. Middle row: plots of the noisy data and of the reconstruction of $p(x)$ at final iteration for $J=25,25\cdot2^9$. Bottom row: plots of the reconstruction of the control $u(x)$ at final iteration for $J=25,25\cdot2^9$ and behavior of the relative error $\frac{\|\overline{\vec{u}}-\vec{u}\|_2^2}{\|\overline{\vec{u}}\|_2^2}$.\label{fig:test2}}
\end{figure}

Let us consider $u(x) = \sin(8x)$, $\forall\,x\in[0,\pi]$, and a fixed value of the noise level $\gamma=0.01$. We show that the method provides a good performance also cases where the control function has a high-frequency profile. In Figure~\ref{fig:test2} we consider the results obtained with $J=25$ and $J=25\cdot2^9=12800$ sampling from the initial distribution $f_0(\vec{u})$. In order to measure the quality of the solution to the inverse problem, we again compare the residual $r$ and the projected residual $R$ (top left plot) and the misfit $\vartheta$ (top right plot) for the two values of $J$. The misfit reaches the noise level in a very small number of iterations for both $J$'s but the residual $r$ for $J=12800$ is reaching a smaller value than the residual computed with $J=25$. This result is observable also in the middle right plot and in the bottom left plot where we compare the reconstruction of $p(x)$ and of the exact control $u(x)$ at the final iteration with the two values of $J$ and using the mean as estimator of the solution. It is very clear that the case with $J=12800$ is providing a better resolution. Finally, in the bottom right plot we show the relative error $\frac{\|\overline{\vec{u}}-\vec{u}\|_2^2}{\|\overline{\vec{u}}\|_2^2}$ as function of the increasing value of $J$ noticing a decreasing behavior.

\subsection{Nonlinear elliptic problem} \label{sec:nonlinear}

\begin{figure}[t!]
	\centering
	\includegraphics[width=0.49\textwidth]{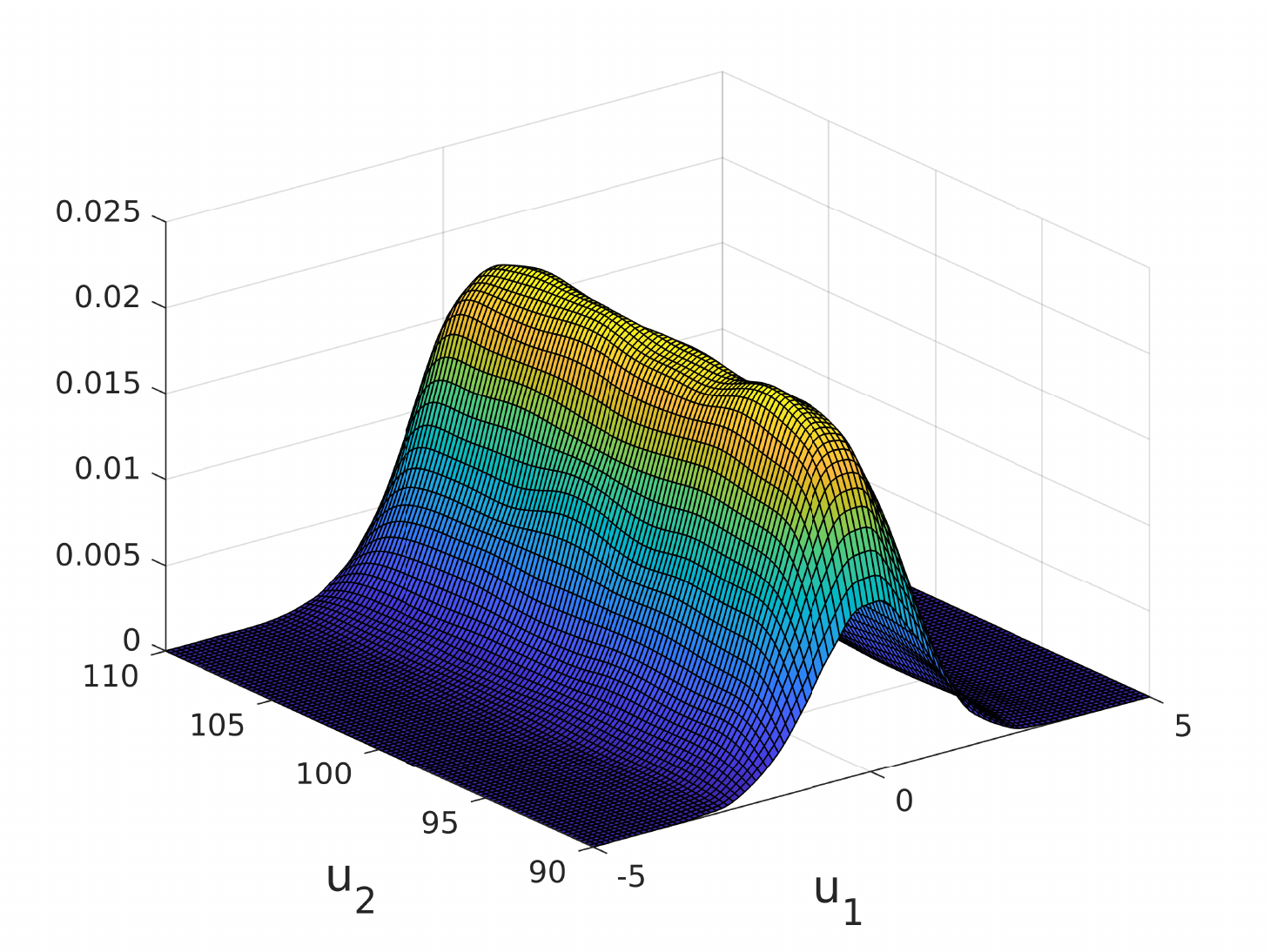}
	\includegraphics[width=0.49\textwidth]{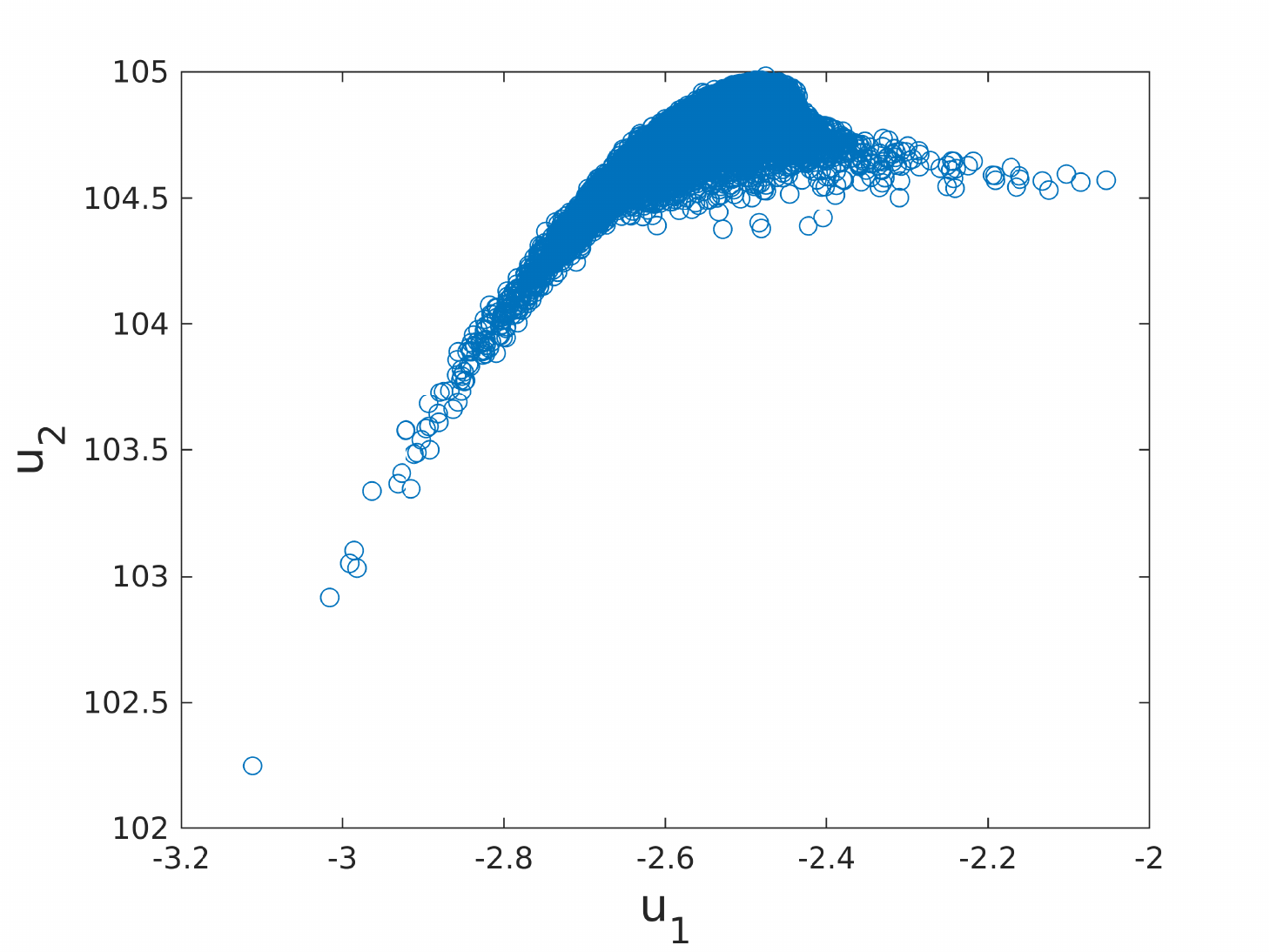}
	\\
	\includegraphics[width=0.49\textwidth]{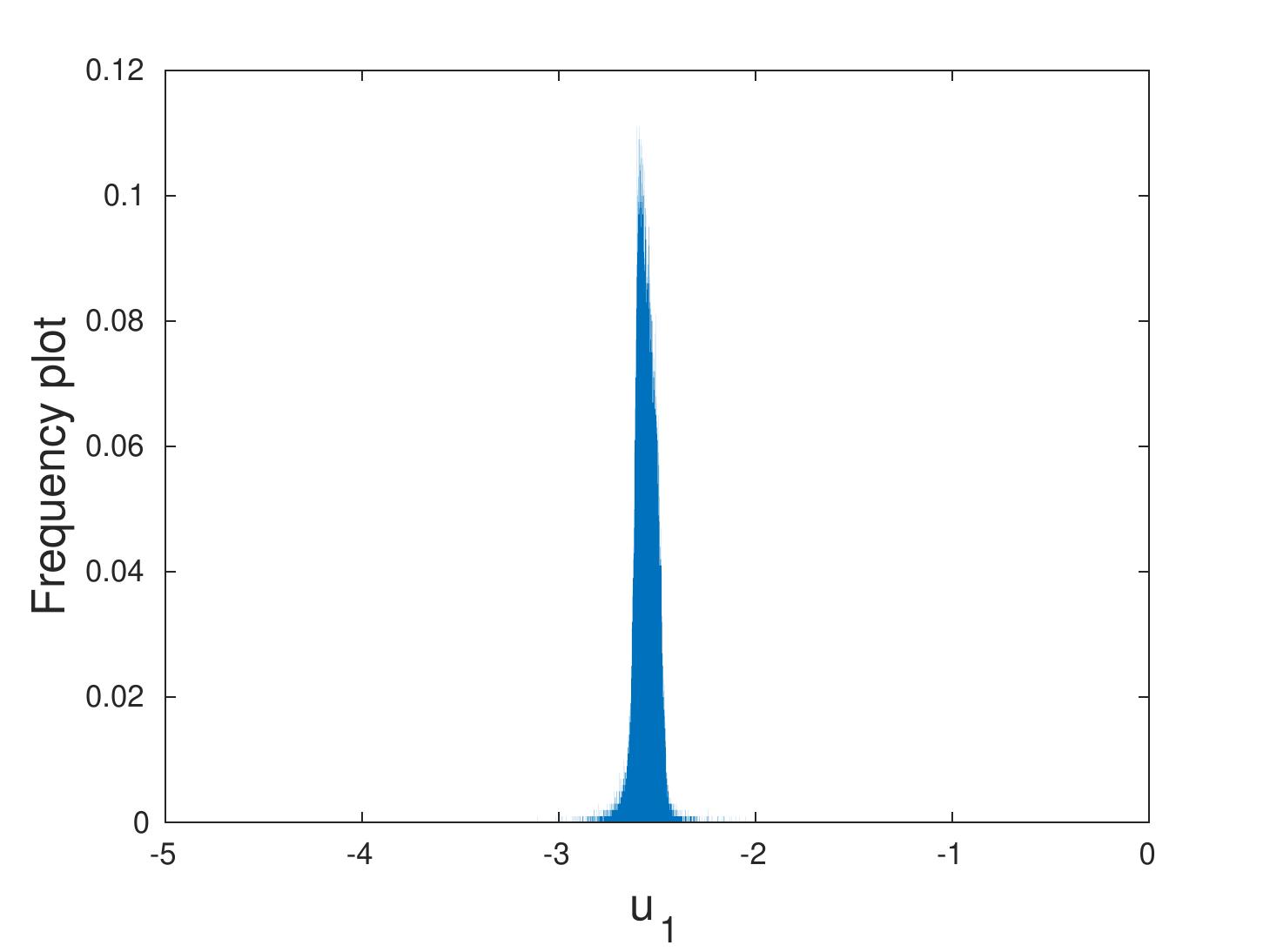}
	\includegraphics[width=0.49\textwidth]{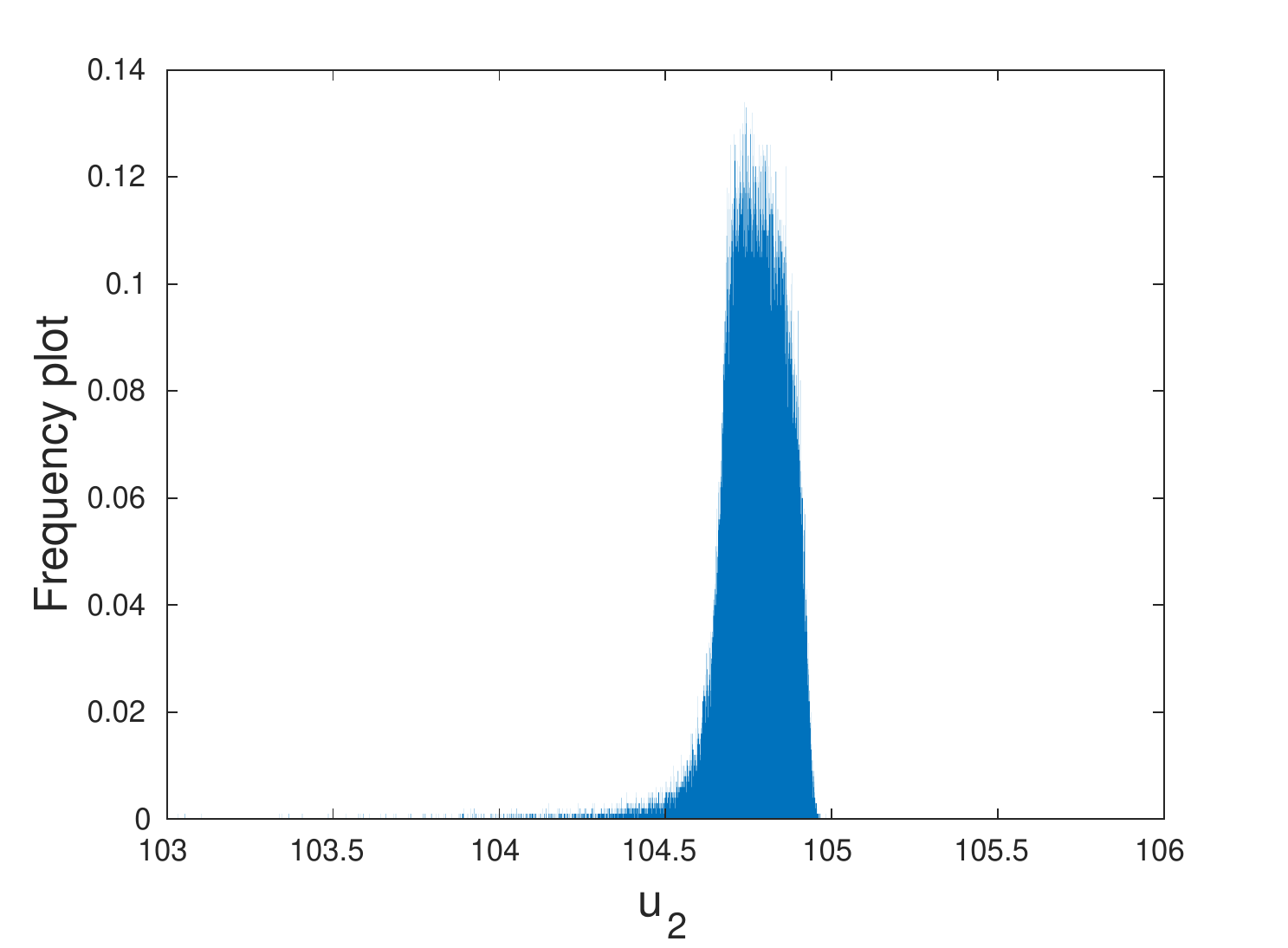}
	\\
	\includegraphics[width=0.49\textwidth]{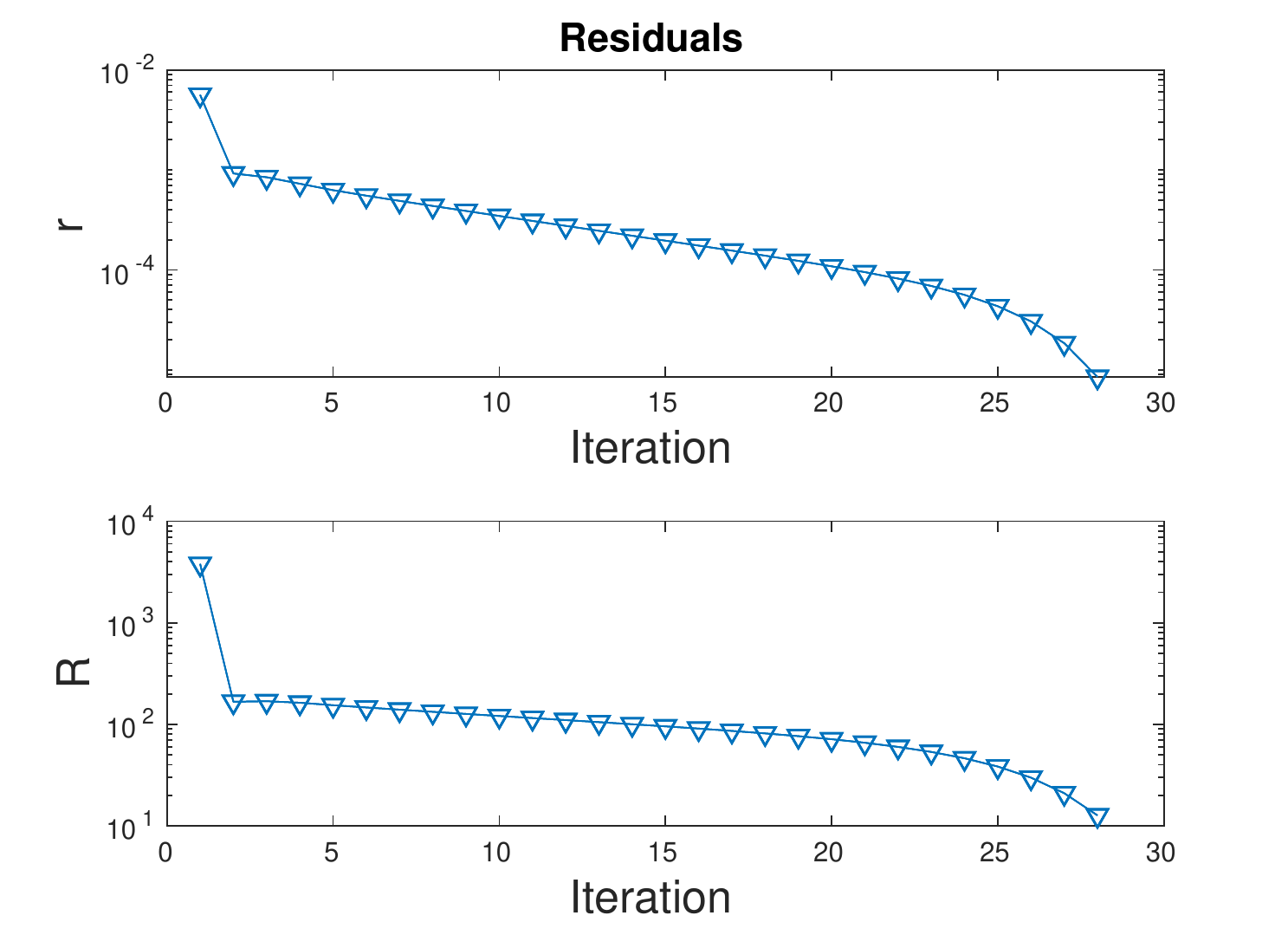}
	\includegraphics[width=0.49\textwidth]{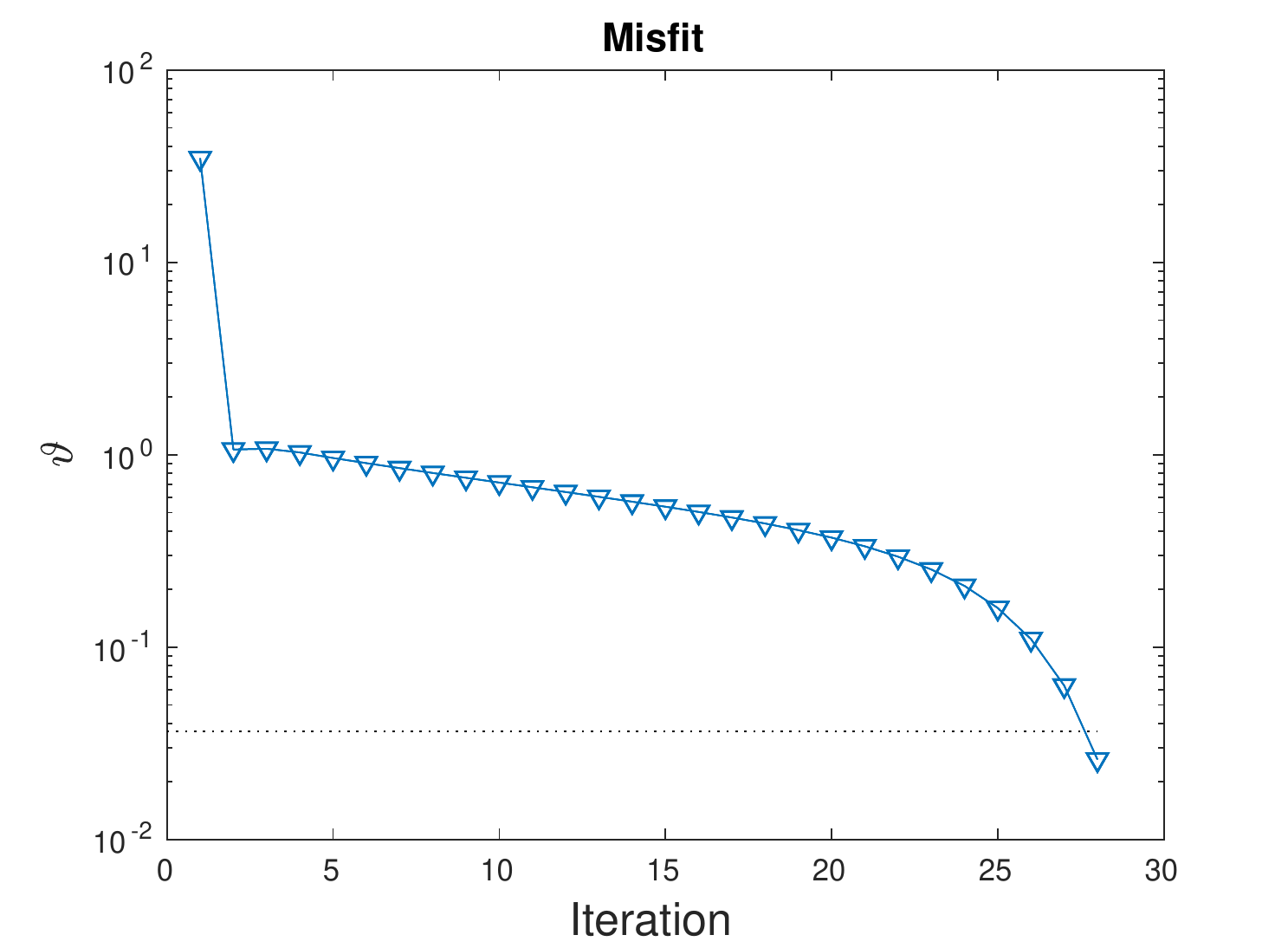}
	\caption{Nonlinear problem. Top row: plots of the density estimation of the initial samples (left) and position of the samples at final iteration (right). Middle row: Marginals of $u_1$ (left) and $u_2$ (right) as relative frequency plot. Bottom row: residual errors $r$ and $R$ (left) and misfit error (right).\label{fig:nonlinearExample}}
\end{figure}

The second numerical experiment is a slightly  modified  example proposed in~\cite{ernstetal2015}. We consider a one-dimensional elliptic boundary value problem given by 
$$
	-\frac{\mathrm{d}}{\mathrm{d}x} \left( \exp(u_1) \frac{\mathrm{d}}{\mathrm{d}x} p(x) \right) = f(x), \quad x\in[0,1]
$$
with boundary conditions $p(0)=p_0$ and $p(1)=u_2$, where $\vec{u}=(u_1,u_2)$ is the unknown control. The exact solution of this problem is given by
$$
	p(x) = p_0 + (u_2-p_0) + \exp(-u_1) \left( -S_x(F) + S_1(F)x \right) 
$$
where $S_x(g)=\int_0^x g(y)\mathrm{d}y$ and $F(x)=S_x(f)=\int_0^x f(y)\mathrm{d}y$. In the following example we consider $f(x)=1$, $\forall\,x\in[0,1]$ and $p_0=0$, so that the explicit solution is given by 
$$
	p(x) = u_2 x + \exp(-u_1) \left( -\frac{x^2}{2} + \frac{x}{2} \right).
$$
We assume to have noisy measurements of $p$ at the points $x_1 = \frac14$ and $x_2=\frac34$ with value $\vec{y} = (27.5,79.7)$. The goal is to seek the control $\vec{u}$ based on the knowledge of $\vec{y}$, of the prior $f_0(\vec{u})$ and of the noise model. More precisely, we consider a prior information given by $\vec{u} \sim \mathcal{N}(0,1) \otimes \mathcal{U}(90,110)$ and a Gaussian white noise $\vecsym{\eta} \sim \mathcal{N}(\vec{0},\gamma^2 \vec{I})$, with $\gamma = 0.1$ and $\vec{I}\in\R^{2\times 2}$ begin the identity matrix. Thus, as in Section~\ref{sec:elliptic}, noisy observations are simulated by
$$
	\vec{y} = \vec{p} + \vecsym{\eta} = \G(\vec{u}^\dagger) + \vecsym{\eta}
$$
where the forward model is defined as
\begin{align*}
	\G \colon \vec{u} \in \R^2 \mapsto \vec{p}=(p(x_1),p(x_2))\in\R^2.
\end{align*}
The example has $d=2$ dimension of the control in order to make a comparison between the solution to the inverse problem provided by the kinetic method and by the Bayes' formula. In particular, it is possible to analyze the approximation of the mean estimator and of the posterior distribution computed by the kinetic equation~\eqref{eq:kineticFP}. However, observe that~\eqref{eq:kineticFP} is derived by assuming a linear forward operator $G$ but in this example the model $\G$ is nonlinear. Thus, inspired by~\eqref{eq:continuousEnKF1}, we consider a small modification of the microscopic interaction rule~\eqref{eq:microInteraction} with $\vec{K} = \vec{I}$ identity matrix given by
$$
\vec{u} = \vec{u}_* + \epsilon \, \vec{C}(\vec{U}_*) \vecsym{\Gamma} (\vec{y}-\G(\vec{u}_*)) + \sqrt{\epsilon} \, \vecsym{\xi}
$$
in order to perform the simulations for the nonlinear model. 

For this example, the true posterior mean is computed in~\cite{ernstetal2015} thanks to  Bayes' formula and it is given by $(-2.65,104.5)$. In Figure~\ref{fig:nonlinearExample} we show the results provided by the kinetic model. The top row shows the density estimation of the $J=10^5$ sampling from the initial distribution (left plot) and the positions of the samples at the last iterations. Again, we use the discrepancy principle as stopping criterion. The middle row shows the marginals of $u_1$ and $u_2$ as relative frequency plots. The solution computed as the mean estimator of the kinetic distribution is $\vec{u} = (-2.56,
104.77)$ is very close to the true posterior mean, as proved also by the plot of the residuals in the bottom left panel of Figure~\ref{fig:nonlinearExample}. The application of the original EnKF method provides $\vec{u} = (-2.92,
105.14)$ which is less accurate, see also~\cite{ernstetal2015}.

\section{Conclusions} \label{sec:conclusion}

In this paper we have introduced a kinetic model for the solution to inverse problems. The kinetic equation has been derived as mean-field limit of the Ensemble Kalman Filter method for infinitely large ensemble. The introduction of a continuous equation describing the evolution of the probability distribution of the unknown control guarantees several advantages: information on statistical quantities of the solution, implicit regularization modeled by the initial distribution, analysis of the properties of the solution.

	The derivation of the kinetic equation has also the advantage to provide a different interpretation of the method and 
	a possibly different scheme using binary collisions with consequent computational gain for numerical simulations. This leads to a different scheme as well as a modified scheme 
	as introduced in the paper. A linear stability analysis for the simple setting of a one dimensional control has showed that the modified method has only stable solutions. Numerical simulations have been performed in order to investigate the good performance of the kinetic equation in providing solutions to inverse problems.


\section*{Acknowledgments}
The authors would like to thank the German Research Foundation DFG for the kind support within the Cluster of Excellence Internet of Production (IoP).\smallskip\\
The authors also acknowledge support by DFG HE5386/14,15.\smallskip\\
Giuseppe Visconti is member of the ``National Group for Scientific Computation (GNCS-INDAM)''.

\bibliographystyle{plain}
\bibliography{references}

%
%

\end{document}